\documentclass[10pt]{amsart}

\usepackage{amsmath,amsfonts,amssymb,amsthm,epsfig,young}
\usepackage[vcentermath]{youngtab}

\def\C{\mathbb{C}}
\def\Z{\mathbb{Z}}
\def\N{\mathbf{N}}

\def\g{\ensuremath{\mathfrak{g}}}

\def\V{\mathbf{V}}
\def\W{\mathbf{W}}
\def\v{\mathbf{v}}
\def\w{\mathbf{w}}

\def\B{\mathbf{B}}
\def\isgn{i\text{-sgn}}
\def\ke{{\tilde e}}
\def\kf{{\tilde f}}

\newcommand{\yl}{14pt}
\newcommand{\yh}{7pt}

\newcommand{\ffbox}[1]{                                        %full box entry
\setbox9=\hbox{$\scriptstyle\overline{1}$}
\framebox[\yl][c]{\rule{0mm}{\ht9}${\scriptstyle #1}$} }

\newcommand{\fhbox}[1]{                                        %half box entry
\setbox9=\hbox{$\scriptstyle\overline{1}$}
\framebox[\yh][c]{\rule{0mm}{\ht9}${\scriptstyle #1}$} }

\newcommand{\fsbox}[1]{                                        %half box entry
\setbox9=\hbox{$\scriptstyle\overline{1}$}
\framebox[\yh][c]{\rule{0mm}{\ht9}${\scriptscriptstyle #1}$} }

\DeclareMathOperator{\Hom}{Hom}

\DeclareMathOperator{\End}{End}
\DeclareMathOperator{\Aut}{Aut}
\DeclareMathOperator{\inc}{in}
\DeclareMathOperator{\out}{out}
\DeclareMathOperator{\tr}{tr}
\DeclareMathOperator{\Irr}{Irr}

\DeclareMathOperator{\wt}{wt}

\DeclareMathOperator{\Coker}{Coker}

\newtheorem{theo}{Theorem}[section]
\newtheorem{prop}[theo]{Proposition}
\newtheorem{lem}[theo]{Lemma}
\newtheorem{cor}[theo]{Corollary}

\newtheorem{defin}[theo]{Definition}
\newtheorem*{rem*}{Remark}
\newtheorem{rem}[theo]{Remark}
\numberwithin{equation}{section}

\begin{document}
\title{Geometric and Combinatorial Realizations of Crystal Graphs}
\author{Alistair Savage}
\address{Max-Planck-Institut f\" ur Mathematik \\ P.O. Box 7280 \\
D-53072 Bonn \\ Germany}
\email{alistair.savage@aya.yale.edu}
\thanks{This research was supported in part by the Natural
Sciences and Engineering Research Council (NSERC) of Canada and
was partially conducted by the author for the Clay Mathematics
Institute.}
\date{October 20, 2003}

\begin{abstract}
For irreducible integrable highest weight modules of the finite
and affine Lie algebras of type $A$ and $D$, we define an
isomorphism between the geometric realization of the crystal
graphs in terms of irreducible components of Nakajima quiver
varieties and the combinatorial realizations in terms of Young
tableaux and Young walls. For type $A_n^{(1)}$, we extend the
Young wall construction to arbitrary level, describing a
combinatorial realization of the crystals in terms of new objects
which we call Young pyramids.
\end{abstract}

\maketitle

\section*{Introduction}
In \cite{FS03}, the author and I. B. Frenkel drew a connection
between two different approaches to representations of affine Lie
algebras of type $A$. This yielded an explicit enumeration of the
irreducible components of certain Nakajima quiver varieties in
terms of Young and Maya diagrams while at the same time yielding
an alternative and much simpler geometric proof of the main result
of \cite{D89} on the construction of bases of representations of
affine Lie algebras.  The geometric proof involved arguments using
commutative diagrams and the Young and Maya diagrams involved in
the enumerations where realized in terms of these.

In \cite{S03}, the author applied the techniques of \cite{FS03} to
the case of the spin representations of type $D$.  Again, an
explicit enumeration of the irreducible components in terms of
Young diagrams was given.  In \cite{S03} and in the simplest cases
appearing in \cite{FS03}, the geometric action was explicitly
computed.  Furthermore, in \cite{S03} the geometric action of the
entire Clifford algebra used in the classical construction of the
spin representations was defined and explicitly computed.

In \cite{KS97,S02}, Kashiwara and Saito endowed the irreducible
components of Lusztig and Nakajima quiver varieties with the
structure of a crystal. In particular, the set of irreducible
components of Nakajima's quiver varieties were given the structure
of the crystals of highest weight irreducible representations.
These crystals have also been realized in purely combinatorial
ways. For instance, for the classical Lie algebras, the crystals
have been realized on tableaux \cite{KN94} and for basic
representations of affine Lie algebras, the crystals have been
realized on Young walls \cite{K03}.  See \cite{HK} for a review of
both of these realizations as well as a review of the theory of
crystals. It is natural to ask what connection it is possible to
make between the geometric and combinatorial realizations of
crystals.  This is the subject of the current paper. For finite
type $A$, we are able to reinterpret the results of \cite{FS03} in
terms of Young tableaux. Specifically, we enumerate the
irreducible components of Nakajima's quiver varieties by the Young
tableaux appearing in the combinatorial crystal realization.
Furthermore, we show this correspondence is actually a crystal
isomorphism.  This approach has the added advantage of being able
to be extended to type $D$ where we also enumerate the irreducible
components by the tableaux appearing in the combinatorial crystal
realization and show that this correspondence is a crystal
isomorphism.

We also consider the affine case.  For both type $A$ and $D$, we
obtain an enumeration of the irreducible components of Nakajima's
quiver varieties for the basic representation by Young walls. For
type $A$ we are able to extend the combinatorial construction
beyond level one to arbitrary level.  We thus define a new
combinatorial object which we call \emph{Young pyramids}, inspired
by the geometry of quiver varieties, realizing the crystals of
arbitrary integrable irreducible highest weight representations of
type $A_n^{(1)}$.

The connection between the geometry of quiver varieties and the
combinatorial objects is very beautiful and has numerous benefits.
First of all, we obtain a very explicit enumeration of irreducible
components of Nakajima's quiver varieties in terms of classical
combinatorial objects.  Secondly, since Nakajima's theory gives us
a realization of Lie algebra representations in terms of the
homology of or constructible functions on the quiver varieties, we
can translate the action of the Lie algebra to the vector space
spanned by the combinatorial objects.  Thus we obtain the full
(rather than just the crystal) structure of the representation in
terms of these objects.  This action seems to be new and in some
cases is the first realization of the entire structure on such
objects. Finally, the geometry naturally explains many
characteristics of the combinatorial crystal realizations. In
particular, entries in a Young tableaux or columns in a Young wall
correspond to specific representations of the quiver. Possible
maps between these representations determine the ordering on the
index set of the tableaux.  Further notions of Young walls such as
the ground state wall and the pattern for building the walls are
also naturally explained by the geometry of quiver varieties.  The
geometric interpretation of classical objects such as Young
tableaux also suggests further avenues of research.  For instance,
it would be interesting to interpret such results as the
Littlewood-Richardson rule in terms of the geometry of
Malkin-Nakajima tensor product varieties.

For other connections between quiver varieties and Young tableaux,
we refer the reader to \cite{N94b,N03} where Nakajima noted a
relationship between the homology of quiver varieties and Young
tableaux for types $A$ and $D$ using methods different than those
employed in the current paper.

The organization of this paper is as follows.  In
Sections~\ref{sec:lus_def} and \ref{sec:def_nak} we review the
definition of Lusztig and Nakajima quiver varieties.  In
Sections~\ref{sec:ca_qv} and \ref{sec:ca_tab} we recall the
realization of various crystal graphs by the geometry of quiver
varieties and the combinatorics of Young tableaux.  We reinterpret
some results of \cite{FS03} on the enumeration of irreducible
components of Nakajima quiver varieties in terms of Young tableaux
in Section~\ref{sec:A_enum} and in Section~\ref{sec:A_ident} we
show that this enumeration actually yields a crystal isomorphism.
In Section~\ref{sec:D_ident} we define the crystal isomorphism
between the geometric and Young tableaux realizations of type $D$
crystals. We introduce the Young pyramid realization of the
crystal graphs of arbitrary irreducible integrable highest weight
representations of type $A_n^{(1)}$ in
Section~\ref{sec:Ahat_ident} and define a crystal isomorphism with
the geometric realization.  In Section~\ref{sec:Dhat_ident} we
define an isomorphism between the geometric and Young wall
realizations of the crystal graphs of the basic representations of
type $D_n^{(1)}$.  Finally, in Section~\ref{sec:paths} we note the
connection to the path model of \cite{D89b,D89} and describe how
our results yield a new action in the space of paths and Young
walls/pyramids.

The author would like to thank I. B. Frenkel for useful
discussions and suggestions and the Max-Planck-Institut f\" ur
Mathematik where much of this research took place.

%%%%%%%%%%%%%%%%%%%%%%%%%%%%%%%%%%%%%%%%%%%%%%%%%%%%%%%%%%%%%%%%%%%%%%

\section{Lusztig's quiver variety}
\label{sec:lus_def}

In this section, we will recount the description given in
\cite{L91} of Lusztig's quiver variety and its irreducible
components. See this reference for details, including proofs.

\subsection{General definitions}

Let $I$ be the set of vertices of the Dynkin graph of a symmetric
Kac-Moody Lie algebra $\mathfrak{g}$ and let $H$ be the set of
pairs consisting of an edge together with an orientation of it.
For $h \in H$, let $\inc(h)$ (resp. $\out(h)$) be the incoming
(resp. outgoing) vertex of $h$.  We define the involution $\bar{\
}: H \to H$ to be the function which takes $h \in H$ to the
element of $H$ consisting of the same edge with opposite
orientation.  An \emph{orientation} of our graph is a choice of a
subset $\Omega \subset H$ such that $\Omega \cup \bar{\Omega} = H$
and $\Omega \cap \bar{\Omega} = \emptyset$.

Let $\mathcal{V}$ be the category of
finite-dimensional $I$-graded vector spaces $\V = \oplus_{i
  \in I} \V_i$ over $\C$ with morphisms being linear maps
respecting the grading.  Then $\V \in \mathcal{V}$ shall denote
that $\V$ is an object of $\mathcal{V}$.  We identify the graded
dimension $\v$ of $\V$ with the element $\sum_{i \in I} \v_i
\alpha_i$ of the root lattice of $\mathfrak{g}$. Here the
$\alpha_i$ are the simple roots corresponding to the vertices of
our quiver (graph with orientation), whose underlying graph is the
Dynkin graph of $\mathfrak{g}$.

Given $\V \in \mathcal{V}$, let
\[
\mathbf{E_V} = \bigoplus_{h \in H} \Hom (\V_{\out(h)},
\V_{\inc(h)}).
\]
For any subset $H' \subset H$, let $\mathbf{E}_{\V, H'}$ be the
subspace of $\mathbf{E_V}$ consisting of all vectors $x = (x_h)$
such that $x_h=0$ whenever $h \not\in H'$.  The algebraic group
$G_\V = \prod_i \Aut(\V_i)$ acts on $\mathbf{E_V}$ and
$\mathbf{E}_{\V, H'}$ by
\[
(g,x) = ((g_i), (x_h)) \mapsto (g_{\inc(h)} x_h g_{\out(h)}^{-1}).
\]

Define the function $\varepsilon : H \to \{-1,1\}$ by $\varepsilon
(h) = 1$ for all $h \in \Omega$ and $\varepsilon(h) = -1$ for all
$h \in {\bar{\Omega}}$.  The Lie algebra of $G_\V$ is
$\mathbf{gl_V} = \prod_i \End(\V_i)$ and it acts on $\mathbf{E_V}$
by
\[
(a,x) = ((a_i), (x_h)) \mapsto [a,x] = (x'_h) = (a_{\inc(h)}x_h - x_h
a_{\out(h)}).
\]
Let $\left<\cdot,\cdot\right>$ be the nondegenerate,
$G_\V$-invariant, symplectic form on
$\mathbf{E_V}$ with values in $\C$ defined by
\[
\left<x,y\right> = \sum_{h \in H} \varepsilon(h) \tr (x_h y_{\bar{h}}).
\]
Note that $\mathbf{E_V}$ can be considered as the cotangent space of
$\mathbf{E}_{\V, \Omega}$ under this form.

The moment map associated to the $G_{\V}$-action on the
symplectic vector space $\mathbf{E_V}$ is the map $\psi : \mathbf{E_V}
\to \mathbf{gl_V}$ with $i$-component $\psi_i : \mathbf{E_V} \to \End
\V_i$ given by
\[
\psi_i(x) = \sum_{h \in H,\, \inc(h)=i} \varepsilon(h) x_h x_{\bar{h}} .
\]

\begin{defin}[\cite{L91}]
\label{def:nilpotent}
An element $x \in \mathbf{E_V}$ is said to be \emph{nilpotent} if
there exists an $N \ge 1$ such that for any sequence $h'_1, h'_2,
\dots, h'_N$ in $H$ satisfying $\out (h'_1) = \inc (h'_2)$, $\out (h'_2) =
\inc (h'_3)$, \dots, $\out (h'_{N-1}) = \inc (h'_N)$, the composition
$x_{h'_1} x_{h'_2} \dots x_{h'_N} : \V_{\out (h'_N)} \to
  \V_{\inc (h'_1)}$ is zero.
\end{defin}

\begin{defin}[\cite{L91}] Let $\Lambda_\V$ be the set of all
  nilpotent elements $x \in \mathbf{E_V}$ such that $\psi_i(x) = 0$
  for all $i \in I$.
\end{defin}

\begin{prop}[\cite{L91}]
\label{prop:irrcomp:Ainfty} For \g\ a symmetric Lie algebra of
finite type, the irreducible components of $\Lambda_\V$ are the
closures of the conormal bundles of the various $G_\V$-orbits in
$\mathbf{E}_{\V, \Omega}$.
\end{prop}

%%%%%%%%%%%%%%%%%%%%%%%%%%%%%%%%%%%%%%%%%%%%%%%%%%%%%%%%%%%%%%%%%%%%%%%

\subsection{Type $A$}
\label{sec:areps}
\begin{figure}
\centering \epsfig{file=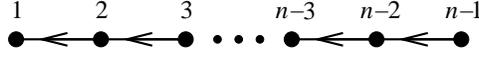,width=0.5\textwidth}
\caption{The quiver of type $A_{n-1}$. \label{fig:quiver_an}}
\end{figure}

Let $\mathfrak{g}=\mathfrak{sl}_n$ be the simple Lie algebra of type
$A_{n-1}$. Let
$I=\{1,\dots,n-1\}$ be the set of vertices of a graph with the set of oriented
edges given by
\begin{gather*}
H=\{h_{i,j} \ |\ i,j \in I,\ |i-j|=1\}.
\end{gather*}
For two adjacent vertices $i$ and $j$, $h_{i,j}$ is the oriented
edge from vertex $i$ to vertex $j$.  Thus $\out(h_{i,j}) = i$ and
$\inc(h_{i,j})=j$. We define the involution $\bar{\ } : H \to H$
as the function that interchanges $h_{i,j}$ and $h_{j,i}$. Let
$\Omega = \{h_{i,i-1}\ |\ 2 \le i \le n-1 \}$.  We picture this
quiver as in Figure~\ref{fig:quiver_an}.

\begin{figure}
\centering \epsfig{file=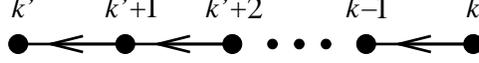,width=0.5\textwidth}
\caption{The string depicting the representation
$(\V(k',k),x(k',k))$. \label{fig:string_an}}
\end{figure}

For two integers $k',k$ such that $1 \le k' \le k \le n-1$, define
$\V(k',k) \in \mathcal{V}$ to be the vector space with basis $\{
e_r\ |\ k' \le r \le k\}$.  We require that $e_r$ has degree $r
\in I$. Let $x(k', k) \in \mathbf{E}_{\V(k',k), \Omega}$ be
defined by $x(k',k) : e_r \mapsto e_{r-1}$ for $k' \le r \le k$,
where $e_{k'-1} = 0$.  It is clear that $(\V(k',k), x(k',k))$ is
an indecomposable representation of our quiver (i.e. element of
$\mathbf{E}_{\V,\Omega}$). Conversely, any indecomposable
finite-dimensional representation $(\V,x)$ of our quiver is
isomorphic to some $(\V(k',k),x(k',k))$. We picture such an
indecomposable representation as a string with endpoints $k'$ and
$k$ (see Figure~\ref{fig:string_an}).

Let $Z$ be the set of all pairs $(k',k)$ of integers such that $1 \le
k' \le k \le n-1$ and
let $\tilde Z$ be
the set of all functions $Z \to \N$ with finite support.

It is easy to see that for $\V \in \mathcal{V}$, the set of
$G_\V$-orbits in
$\mathbf{E}_{\V, \Omega}$ is naturally indexed by the subset
$\tilde Z_\V$ of $\tilde Z$ consisting of those
$f \in \tilde Z$ such that
\[
\sum_{k' \le i \le k} f(k',k) = \dim \V_i
\]
for all $i \in I$.  Here the sum is over all $k', k$ such that $1
\le k' \le i \le k \le n-1$.  Corresponding to a given $f$ is the
orbit consisting of all representations isomorphic to a sum of the
indecomposable representations $(V(k',k),x(k',k))$, each occuring
with multiplicity $f(k',k)$. Denote by $\mathcal{O}_f$ the
$G_\V$-orbit corresponding to $f \in \tilde Z_\V$.  Let
$\mathcal{C}_f$ be the conormal bundle to $\mathcal{O}_f$ and let
$\bar{\mathcal{C}}_f$ be its closure.  We then have the following
proposition.

\begin{prop}
The map $f \to \bar{\mathcal{C}}_f$ is a one-to-one correspondence
between the set
$\tilde Z_\V$ and the set of irreducible
components of $\Lambda_\V$.
\end{prop}
\begin{proof}
This follows immediately from Proposition~\ref{prop:irrcomp:Ainfty}.
\end{proof}

%%%%%%%%%%%%%%%%%%%%%%%%%%%%%%%%%%%%%%%%%%%%%%%%%%%%%%%%%%%%%%%%%%%%%%%

\subsection{Type $D$}

\begin{figure}
\centering \epsfig{file=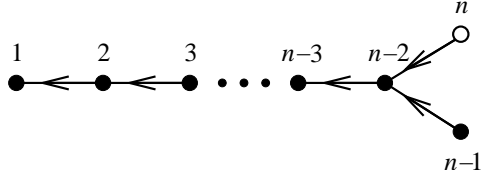,width=0.5\textwidth}
\caption{The quiver of type $D_n$. We represent the $n$th vertex
by an open dot to distinguish it from the $(n-1)$st vertex.
\label{fig:quiver_dn}}
\end{figure}

Now let $\g=\mathfrak{so}_{2n}\C$ be the simple Lie algebra of
type $D_n$. Let $I=\{1,2,\dots,n\}$ be the set of vertices of the
Dynkin graph of \g, labelled and oriented as in
Figure~\ref{fig:quiver_dn}.
As for type $A$, we label each oriented edge by its
incoming and outgoing vertices. That is, if vertices $i$ and $j$
are connected by an edge, $h_{i,j}$ denotes the oriented edge with
$\out(h)=i$, $\inc(h)=j$.

%%%%%%%%%%%%%%%%%%%%%%%%%%%%%%%%%%%%%%%%%%%%%%%%%%%%%%%%%%%%%%%%%%%%%%%

\section{Nakajima's quiver variety}
\label{sec:def_nak}

We introduce here a description of the quiver varieties first
presented in \cite{N94}.  See \cite{N94} and \cite{N98} for details.

\begin{defin}[\cite{N94}]
\label{def:lambda} For $\v, \w \in (\Z_{\ge 0})^I$, choose
$I$-graded vector spaces $\V$ and $\W$ of graded dimensions $\v$
and $\w$ respectively.  We associate $\w$ with the element $\sum_i
\w_i \omega_i$ of the weight lattice of \g, where the $\omega_i$
are the fundamental weights of \g.  Recall that we identified $\v$
with the weight $\sum_i \v_i \alpha_i$.  Then define
\[
\Lambda \equiv \Lambda(\v,\w) =
\Lambda_\V \times \sum_{i \in I} \Hom (\V_i, \W_i).
\]
\end{defin}

Now, suppose that $\mathbf{S}$ is an $I$-graded subspace of $\V$.
For $x \in
\Lambda_\V$ we say that $\mathbf{S}$ is
\emph{$x$-stable} if $x(\mathbf{S}) \subset \mathbf{S}$.

\begin{defin}[\cite{N94}]
\label{def:lambda-stable} Let $\Lambda^{\text{st}} =
\Lambda(\v,\w)^{\text{st}}$ be the set of all $(x, t) \in
\Lambda(\v,\w)$ satisfying the following condition:  If
$\mathbf{S}=(\mathbf{S}_i)$ with $\mathbf{S}_i \subset \V_i$ is
$x$-stable and $t_i(\mathbf{S}_i) = 0$ for $i \in I$, then
$\mathbf{S}_i = 0$ for $i \in I$.
\end{defin}

The group $G_\V$ acts on $\Lambda(\v,\w)$ via
\[
(g,(x,t)) = ((g_i), ((x_h), (t_i))) \mapsto ((g_{\inc (h)} x_h
g_{\out (h)}^{-1}), (t_i g_i^{-1})).
\]
and the stabilizer of any point of
$\Lambda(\v,\w)^{\text{st}}$ in $G_{\V}$ is trivial
  (see \cite[Lemma~3.10]{N98}).  We then make the following definition.
\begin{defin}[\cite{N94}]
\label{def:L}
Let $\mathcal{L} \equiv \mathcal{L}(\v,\w) =
\Lambda(\v,\w)^{\text{st}} / G_{\V}$.
\end{defin}

Let $\Irr \mathcal{L}(\v,\w)$ (resp. $\Irr
\Lambda(\v,\w)$) be the set of irreducible components of
$\mathcal{L}(\v,\w)$ (resp. $\Lambda(\v,\w)$).
Then $\Irr \mathcal{L}(\v,\w)$ can be identified with
\[
\{ Y \in \Irr \Lambda(\v,\w)\, |\, Y \cap
\Lambda(\v,\w)^{\text{st}} \ne \emptyset \}.
\]
Specifically, the irreducible components of $\Irr
\mathcal{L}(\v,\w)$ are precisely those
\[
X_f \stackrel{\text{def}}{=} \left( \left( \bar{\mathcal{C}}_f \times
    \sum_{i \in I} \Hom (\V_i, \W_i) \right) \cap
\Lambda(\v,\w)^{\text{st}} \right) / G_\V
\]
which are nonempty.

The following will be used in the sequel.
\begin{lem}
\label{lem:irrcomp}
One has
\[
\Lambda^{\text{st}} = \left\{ (x,t) \in \Lambda\, \left| \, \left(
\bigcap_{h\, :\, \out(h) =i} \ker x_h \right) \cap \ker t_i = 0 \
\forall i \right. \right\}.
\]
\end{lem}
\begin{proof}
Since each $\left( \bigcap_{h\, :\, \out(h) =i} \ker x_h \right)
\cap \ker t_i$ is $x$-stable, the lefthand side is obviously
contained in the righthand side. Now suppose $x$ is an element of
the righthand side.  Let $\mathbf{S} = (\mathbf{S}_i)$ with
$\mathbf{S}_i \subset \V_i$ be $x$-stable and $t_i (\mathbf{S}_i)
= 0$ for $i \in I$.  Assume that $\mathbf{S} \ne 0$. Since all
elements of $\Lambda$ are nilpotent, we can find a minimal value
of $N$ such that the condition in Definition~\ref{def:nilpotent}
is satisfied. Then we can find a $v \in \mathbf{S}_i$ for some $i$
and a sequence $h_1', h_2', \dots, h_{N-1}'$ (empty if $N=1$) in
$H$ such that $\out (h_1') = \inc (h_2')$, $\out (h_2') = \inc
(h_3')$, \dots, $\out (h_{N-2}') = \inc (h_{N-1}')$ and $v' =
x_{h_1'} x_{h_2'} \dots x_{h_{N-1}'} (v) \ne 0$. Now, $v' \in
\mathbf{S}_{i'}$ for some $i' \in I$ by the stability of
$\mathbf{S}$ (hence $t_{i'}(v') = 0$) and $v' \in \bigcap_{h\, :\,
\out(h) =i'} \ker x_h$ by our choice of $N$. This contradicts the
fact that $x$ is an element of the righthand side.
\end{proof}

%%%%%%%%%%%%%%%%%%%%%%%%%%%%%%%%%%%%%%%%%%%%%%%%%%%%%%%%%%%%%%%%%%%%%%%

\section{Crystal action on quiver varieties}
\label{sec:ca_qv} In this section, we review the realization of
the crystal graph of integrable highest weight representations of
a Kac-Moody algebra \g\ with symmetric Cartan matrix via quiver
varieties. See \cite{KS97,S02} for details, including proofs.

Let $\mathbf{w, v, v', v''} \in (\Z_{\ge 0})^I$ be such that $\v =
\mathbf{v'} + \mathbf{v''}$.  Consider the maps
\begin{equation}
\label{eq:diag_action} \Lambda(\v'',\mathbf{0}) \times
\Lambda(\v',\w) \stackrel{p_1}{\leftarrow} \mathbf{\tilde F
(v,w;v'')} \stackrel{p_2}{\rightarrow} \mathbf{F(v,w;v'')}
\stackrel{p_3}{\rightarrow} \Lambda(\v,\w),
\end{equation}
where the notation is as follows.  A point of
$\mathbf{F(v,w;v'')}$ is a point $(x,t) \in \Lambda(\v,\w)$
together with an $I$-graded, $x$-stable subspace $\mathbf{S}$ of
$\V$ such that $\dim \mathbf{S} = \mathbf{v'} = \v - \v''$.  A
point of $\mathbf{\tilde
  F (v,w;v'')}$ is a point $(x,t,\mathbf{S})$ of $\mathbf{F(v,w;v'')}$
together with a collection of isomorphisms $R'_i : \V'_i \cong
\mathbf{S}_i$ and $R''_i : \V''_i \cong \V_i / \mathbf{S}_i$ for
each $i \in I$.  Then we define $p_2(x,t,\mathbf{S}, R',R'') =
(x,t,\mathbf{S})$, $p_3(x,t,\mathbf{S}) = (x,t)$ and
$p_1(x,t,\mathbf{S},R',R'') = (x'',x',t')$ where $x'', x', t'$ are
determined by
\begin{align*}
R'_{\inc(h)} x'_h &= x_h R'_{\out(h)} : \V'_{\out(h)} \to
\mathbf{S}_{\inc(h)}, \\
t'_i &= t_i R'_i : \V'_i \to \W_i \\
R''_{\inc(h)} x''_h &= x_h R''_{\out(h)} : \V''_{\out(h)} \to
\V_{\inc(h)} / \mathbf{S}_{\inc(h)}.
\end{align*}
It follows that $x'$ and $x''$ are nilpotent.

\begin{lem}[{\cite[Lemma 10.3]{N94}}]
One has
\[
(p_3 \circ p_2)^{-1} (\Lambda(\v,\w)^{\text{st}}) \subset p_1^{-1}
(\Lambda(\v'',\mathbf{0}) \times \Lambda(\v',\w)^{\text{st}}).
\]
\end{lem}

Thus, we can restrict \eqref{eq:diag_action} to
$\Lambda^{\text{st}}$, forget the
$\Lambda(\v'',\mathbf{0})$-factor and consider the quotient by
$G_\V$, $G_\mathbf{V'}$.  This yields the diagram
\begin{equation}
\label{eq:diag_action_mod} \mathcal{L}(\v', \w)
\stackrel{\pi_1}{\leftarrow} \mathcal{F}(\v, \w; \v - \mathbf{v'})
\stackrel{\pi_2}{\rightarrow} \mathcal{L}(\v, \w),
\end{equation}
where
\[
\mathcal{F}(\v, \w; \v - \mathbf{v'}) \stackrel{\text{def}}{=} \{
(x,t,\mathbf{S}) \in \mathbf{F(\v,\w;\v-\v')}\,
  |\, (x,t) \in \Lambda(\v,\w)^{\text{st}} \} / G_\V.
\]

For $i \in I$ define $\varepsilon_i : \Lambda(\v, \w) \to \Z_{\ge
0}$ by
\[
\varepsilon_i((x,t)) = \dim_\C \Coker \left( \bigoplus_{h\, :\,
\inc(h)=i} V_{\out(h)} \stackrel{(x_h)}{\longrightarrow} V_i
\right).
\]
Then, for $c \in \Z_{\ge 0}$, let
\[
\mathcal{L}(\v,\w)_{i,c} = \{[x,t] \in \mathcal{L}(\v,\w)\ |\
\varepsilon_i((x,t)) = c\}
\]
where $[x,t]$ denotes the $G_\V$-orbit through the point $(x,t)$.
$\mathcal{L}(\v,\w)_{i,c}$ is a locally closed subvariety of
$\mathcal{L}(\v,\w)$.

Assume $\mathcal{L}(\v,\w)_{i,c} \ne \emptyset$ and let $\v' = \v
- c\mathbf{e}^i$ where $\mathbf{e}^i_j = \delta_{ij}$.  Then
\[
\pi_1^{-1}(\mathcal{L}(\v',\w)_{i,0}) =
\pi_2^{-1}(\mathcal{L}(\v,\w)_{i,c}).
\]
Let
\[
\mathcal{F}(\v,\w;c\mathbf{e}^i)_{i,0} =
\pi_1^{-1}(\mathcal{L}(\v',\w)_{i,0}) =
\pi_2^{-1}(\mathcal{L}(\v,\w)_{i,c}).
\]
We then have the following diagram.
\begin{equation}
\label{eq:crystal-action} \mathcal{L}(\v',\w)_{i,0}
\stackrel{\pi_1}{\longleftarrow} \mathcal{F}(\v,\w;
c\mathbf{e}^i)_{i,0} \stackrel{\pi_2}{\longrightarrow}
\mathcal{L}(\v,\w)_{i,c}
\end{equation}
The restriction of $\pi_2$ to $\mathcal{F}(\v,\w;
c\mathbf{e}^i)_{i,0}$ is an isomorphism since the only possible
choice for the subspace $\mathbf{S}$ of $\V$ is to have
$\mathbf{S}_j = \mathbf{V}_j$ for $j \ne i$ and $\mathbf{S}_i$
equal to the sum of the images of the $x_h$ with $\inc(h)=i$.
$\mathcal{L}(\v',\w)_{i,0}$ is an open subvariety of
$\mathcal{L}(\v',\w)$.

\begin{lem}[\cite{S02}]
\begin{enumerate}
\item For any $i \in I$,
\[
\mathcal{L}(\mathbf{0},\w)_{i,c} =
\begin{cases}
pt & \text{if $c=0$} \\
\emptyset & \text{if $c >0$}
\end{cases}.
\]
\item Suppose $\mathcal{L}(\v,\w)_{i,c} \ne \emptyset$ and $\v' =
\v - c\mathbf{e}^i$.  Then the fiber of the restriction of $\pi_1$
to $\mathcal{F}(\v, \w; c\mathbf{e}^i)_{i,0}$ is isomorphic to a
Grassmanian variety.
\end{enumerate}
\end{lem}

\begin{cor}
\label{cor:irrcomp-isom} Suppose $\mathcal{L}(\v,\w)_{i,c} \ne
\emptyset$.  Then there is a 1-1 correspondence between the set of
irreducible components of $\mathcal{L}(\v - c\mathbf{e}^i,
\w)_{i,0}$ and the set of irreducible components of
$\mathcal{L}(\v, \w)_{i,c}$.
\end{cor}

Let $B(\v,\w)$ denote the set of irreducible components of
$\mathcal{L}(\v,\w)$ and let $B(\w) = \bigsqcup_\v B(\v,\w)$. For
$X \in B(\v,\w)$, let $\varepsilon_i(X) = \varepsilon_i((x,t))$
for a generic point $[x,t] \in X$.  Then for $c \in \Z_{\ge 0}$
define
\[
B(\v,\w)_{i,c} = \{X \in B(\v,\w)\ |\ \varepsilon_i(X) = c\}.
\]
Then by Corollary~\ref{cor:irrcomp-isom}, $B(\v -
c\mathbf{e}^i,\w)_{i,0} \cong B(\v, \w)_{i,c}$.

Suppose that ${\bar X} \in B(\v - c\mathbf{e}^i,\w)_{i,0}$
corresponds to $X \in B(\v,\w)_{i,c}$ by the above isomorphism.
Then we define maps
\begin{gather*}
\kf_i^c : B(\v - c\mathbf{e}^i,\w)_{i,0} \to B(\v,\w)_{i,c},\quad
\kf_i^c({\bar X})
= X, \\
\ke_i^c : B(\v,\w)_{i,c} \to B(\v - c\mathbf{e}^i,\w)_{i,0},\quad
\ke_i^c(X) = {\bar X}.
\end{gather*}
We then define the maps
\[
\ke_i, \kf_i : B(\w) \to B(\w) \sqcup \{0\}
\]
by
\begin{gather*}
\ke_i : B(\v,\w)_{i,c} \stackrel{\ke_i^c}{\longrightarrow} B(\v -
c\mathbf{e}^i, \w)_{i,0} \stackrel{\kf_i^{c-1}}{\longrightarrow}
B(\v -
\mathbf{e}^i, \w)_{i,c-1}, \\
\kf_i : B(\v,\w)_{i,c} \stackrel{\ke_i^c}{\longrightarrow} B(\v -
c\mathbf{e}^i, \w)_{i,0} \stackrel{\kf_i^{c+1}}{\longrightarrow}
B(\v + \mathbf{e}^i, \w)_{i,c+1}.
\end{gather*}
We set $\ke_i(X)=0$ for $X \in B(\v,\w)_{i,0}$ and $\kf_i(X)=0$
for $X \in B(\v,\w)_{i,c}$ with $B(\v,\w)_{i,c+1} = \emptyset$.
We also define
\begin{gather*}
\wt : B(\w) \to P,\quad \wt(X) = \sum_{i\in I} \left( \mathbf{w}_i
\omega_i - \mathbf{v}_i \alpha_i \right) \text{
  for } X \in B(\v,\w), \\
\varphi_i(X) = \varepsilon_i(X) + \left< h_i, \wt(X) \right>.
\end{gather*}

Recall that we can consider $\w$ to be an dominant integral weight
by $\mathbf{w} = \sum_i \w_i \omega_i$.
\begin{prop}[\cite{S02}]
$B(\w)$ is a crystal and is isomorphic to the crystal of the
highest weight $U_q(\g)$-module with highest weight $\w$.
\end{prop}

%%%%%%%%%%%%%%%%%%%%%%%%%%%%%%%%%%%%%%%%%%%%%%%%%%%%%%%%%%%%%%%%%%%%%

\section{Crystal action on tableaux}
\label{sec:ca_tab} We now review the realization of the crystal
graph of finite dimensional representations of Lie algebras of
type $A$ and $D$ via Young tableaux. See \cite{KN94,HK} for
details, including proofs.

Given a set of crystal graphs, the tensor product rule for
crystals gives a very explicit description of the action of the
Kashiwara operators on the multifold tensor product of these
graphs.  Let $\mathcal{B}_j$ be \g-crystals for $i=1, \dots, N$.
Fix $i \in I$ and let $b = b_1 \otimes \dots \otimes b_N \in
\mathcal{B}_1 \otimes \dots \otimes \mathcal{B}_N$.  To each
$b_j$, assign a series of $-$'s and $+$'s with as many $-$'s as
$\varepsilon_i(b_j)$ followed by as many $+$'s as
$\varphi_i(b_j)$.  In the sequence obtained by concatenating the
series for the individual $b_j$'s, cancel all $(+,-)$ pairs to
obtain a sequence $\isgn(b)$ of $-$'s followed by $+$'s. Then the
tensor product rule tells us that $\ke_i$ acts on the $b_j$
corresponding to the rightmost $-$ in $\isgn(b)$ and ${\tilde
f}_i$ acts on the $b_j$ corresponding to the leftmost $+$ in
$\isgn(b)$.

\subsection{Type A}

Let $\g = \mathfrak{sl}_n$ be the Lie algebra of type $A_{n-1}$.
Recall that \g\ acting on the space $\C^n$ by left multiplication
yields the \emph{vector representation}. It has crystal
$\B=\{\ffbox{j}\ |\ j=1,\dots,n\}$ with crystal graph
\[
\ffbox{1} \stackrel{1}{\longrightarrow} \ffbox{2}
\stackrel{2}{\longrightarrow} \dots
\stackrel{n-2}{\longrightarrow} \ffbox{n\!-\!1}
\stackrel{n-1}{\longrightarrow} \ffbox{n} \ .
\]
and $\wt(\ffbox{j})=\epsilon_j$ for $j=1,\dots,n$.  Here the
Cartan subalgebra $\mathfrak{h}$ of $\g$ is the set of traceless
diagonal matrices and $\epsilon_j$ is the element of
$\mathfrak{h}^*$ such that $\epsilon_j(e_{i,i})=\delta_{i,j}$.

A tableau with $m$ boxes represents an element of the tensor
product crystal $\B^{\otimes m}$ by the \emph{Far-Eastern reading}
which proceeds down columns from top to bottom and from right to
left. For example
\[
\young(1123,234,44,5) = \ffbox{3} \otimes \ffbox{2} \otimes
\ffbox{4} \otimes \ffbox{1} \otimes \ffbox{3} \otimes \ffbox{4}
\otimes \ffbox{1} \otimes \ffbox{2} \otimes \ffbox{4} \otimes
\ffbox{5} \ .
\]
From now on, when we say that one box or entry in a tableau is
earlier or later than another, we are using this ordering. Earlier
means ``to the left of" and later means ``to the right of."  We
avoid the words left and right to prevent confusion with the
spatial location of boxes in the tableau.  To denote the spatial
arrangement of boxes or entries in a tableau, we will use the
compass directions north, south, east and west. North is up on the
page, east is right, etc. Now since
\[
\varepsilon_i(\ffbox{i\!+\!1}) = 1,\ \varphi_i(\ffbox{i})=1
\]
and $\varepsilon_i$ and $\varphi_i$ take the value zero on all
other elements of $\B$, the tensor product rule tells us that to
compute the action of the Kashiwara operators $\ke_i$ and $\kf_i$
on a given tableau we neglect all entries not equal to $i$ or
$i+1$ and cancel ($i$,$i+1$) pairs (in the ordering given by the
Far-Eastern reading).  We will say that such a pair of entries are
\emph{$i$-matched}.  Then $\kf_i$ acts by changing the earliest
remaining $i$ to an $i+1$ and $\ke_i$ acts by changing the latest
remaining $i+1$ to an $i$.  If no $i$ (resp. $i+1$) entries
remain, then $\kf_i$ (resp. $\ke_i$) kills the tableau.

As an example, consider the action of $\ke_2$ and
$\kf_2$ on the tableau
\[
T = \young(12233,234,45,5)\ .
\]
Neglecting all entries but 2 and 3, we obtain
\[
\young(\hfill 2233,23\hfill,\hfill \hfill,\hfill)\ .
\]
Read off in the Far-Eastern reading, the entries are $(3,3,2,2,3,2)$.  We
then cancel the 2-matched $(2,3)$ pair to leave the tableau
\[
\young(\hfill \hfill 233,2\hfill \hfill,\hfill \hfill,\hfill)\ .
\]
Thus
\[
\kf_2(T) = \young(12333,234,45,5),\quad \ke_2(T) =
\young(12223,234,45,5).
\]

%%%%%%%%%%%%%%%%%%%%%%%%%%%%%%%%%%%%%%%%%%%%%%%%%%%%%%%%%%%%%%%%%%%%%%

\subsection{Type D}

Let $\g = \mathfrak{so}_{2n}\C$ be the simple Lie algebra of type
$D_n$.  Let $\epsilon_i : M_{2n \times 2n}(\C) \to \C$ be the
linear functional defined by $\epsilon_i(T) = T_{ii}$.  Then the
simple roots and fundamental weights are
\begin{align*}
\alpha_i &= \epsilon_i - \epsilon_{i+1},\quad 1 \le i \le n-1,\\
\alpha_n &= \epsilon_{n-1} + \epsilon_n, \\
\omega_i &= \epsilon_1 + \dots + \epsilon_i,\quad 1 \le i \le n-2,\\
\omega_{n-1} &= \frac{1}{2}(\epsilon_1 + \dots + \epsilon_{n-1} -
\epsilon_n), \\
\omega_n &= \frac{1}{2}(\epsilon_1 + \dots + \epsilon_{n-1} +
\epsilon_n).
\end{align*}

Let $\N = \{1,2,\dots,n,{\bar n},\dots,{\bar 2},{\bar 1}\}$.  We
put a linear ordering on $\N$ by
\[
1 \prec 2 \prec \dots \prec \begin{matrix}n \\ \bar n \end{matrix}
\prec \dots \prec {\bar 2} \prec {\bar 1}
\]
where the order between $n$ and $\bar n$ is not defined. Will will
use roman characters (e.g. $i$ and $j$) to refer to the elements
$1,\dots,n \in \N$ and greek characters (e.g $\alpha$ and $\beta$)
to refer to arbitrary elements of $\N$. Recall that $\g$ acting on
the space $\C^{2n}$ by left multiplication yields the \emph{vector
representation}. It has crystal $\B = \{\ffbox{\alpha}\ |\ \alpha
\in \N\}$ with the following crystal graph.

\begin{center}
\epsfig{file=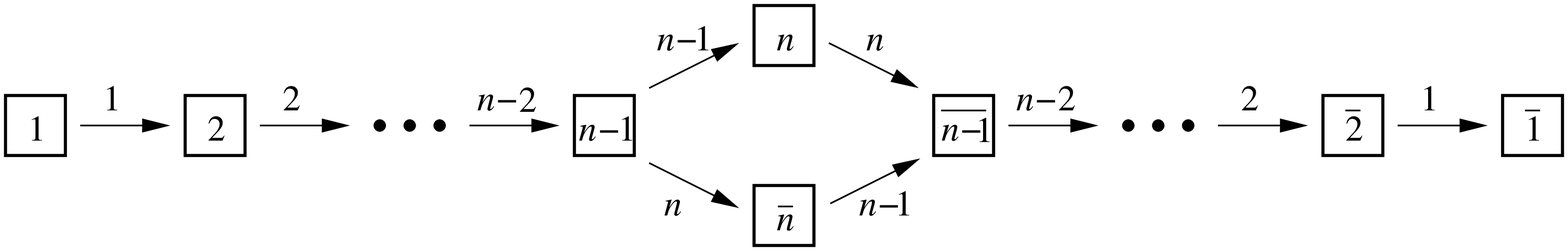,width=\textwidth}
\end{center}

We have that
\[
\wt(\ffbox{j})=\epsilon_j,\quad \wt(\ffbox{\bar j}) = -\epsilon_j
\quad 1 \le j \le n.
\]

Now let $\V_{sp}^+$ and $\V_{sp}^-$ be the spin representations
$V(\omega_n)$ and $V(\omega_{n-1})$ respectively.  The crystal
graphs $\B_{sp}^\pm$ of these representations can be realized
using the (generalized) Young tableaux consisting of half boxes:
\begin{align*}
\B_{sp}^+ &= \left\{ \left.
\raisebox{-0.35in}{\epsfig{file=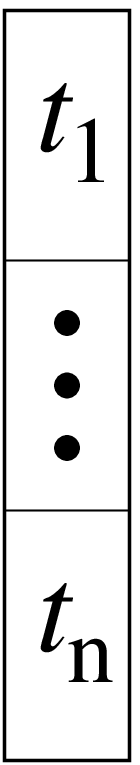, height=0.8in}}\
\right| \
\begin{matrix}
t_i \in \N,\; t_1 \prec \dots \prec t_n, \\
\text{$i$ and $\bar i$ do not appear simultaneously,} \\
\text{if $t_k=n$, then $n-k$ is even,} \\
\text{if $t_k={\bar n}$, then $n-k$ is odd.}
\end{matrix}
\right\},\\
\B_{sp}^- &= \left\{ \left.
\raisebox{-0.35in}{\epsfig{file=halftab_gen.eps, height=0.8in}}\
\right| \
\begin{matrix}
t_i \in \N,\; t_1 \prec \dots \prec t_n, \\
\text{$i$ and $\bar i$ do not appear simultaneously,} \\
\text{if $t_k=n$, then $n-k$ is odd,} \\
\text{if $t_k={\bar n}$, then $n-k$ is even.}
\end{matrix}
\right\}.
\end{align*}

The action of the Kashiwara operators is as follows.
\begin{equation}
\label{eq:kash_halfcol}
\raisebox{-0.7in}{\epsfig{file=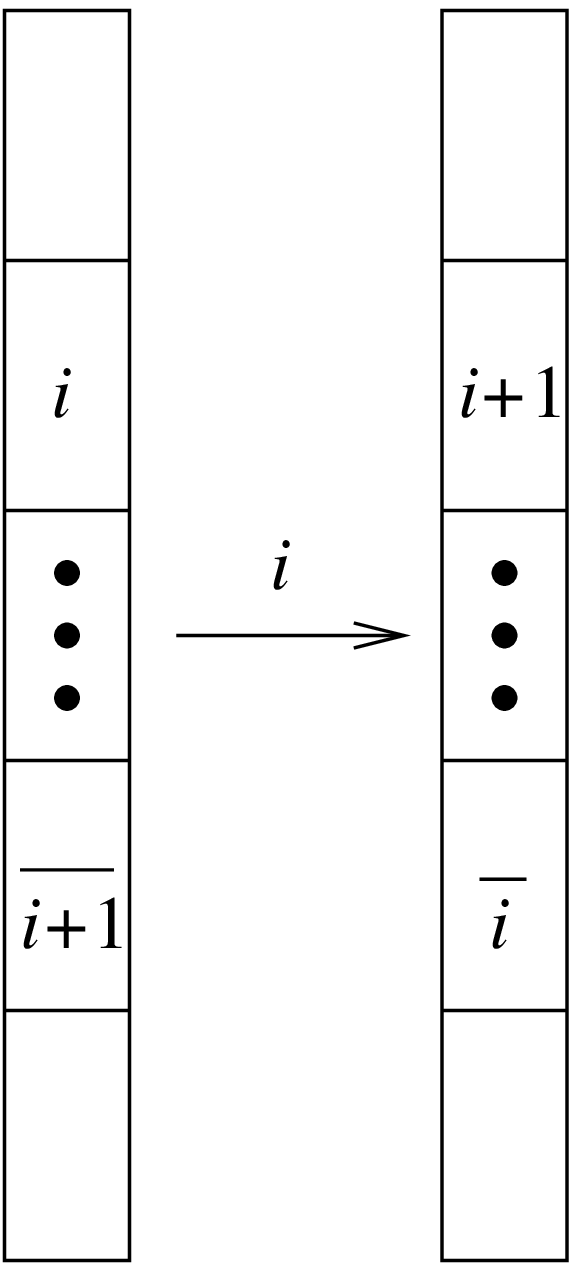, height=1.4in}}
,\ i \ne n,\quad
\raisebox{-0.7in}{\epsfig{file=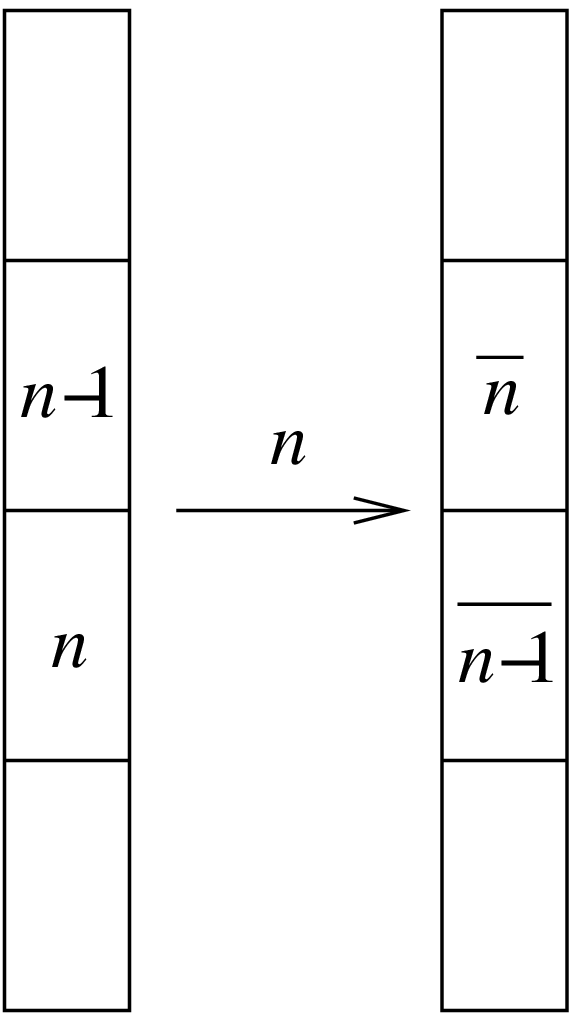,
height=1.4in}}.
\end{equation}

Let $\lambda = \w_1 \omega_1 + \dots + \w_n \omega_n$ ($\w_i \in
\Z_{\ge 0}$) be a dominant integral weight.  Then $\lambda =
\lambda_1 \epsilon_1 + \dots + \lambda_n \epsilon_n$ where
\begin{align*}
\lambda_1 &= \w_1 + \w_2 + \dots + \w_{n-2} + \frac{1}{2}(\w_{n-1}
+ \w_n), \\
\lambda_2 &= \w_2 + \dots + \frac{1}{2}(\w_{n-1}+\w_n), \\
&\ \vdots \\
\lambda_{n-1} &= \frac{1}{2}(\w_{n-1} + \w_n), \\
\lambda_n &= \frac{1}{2}(\w_n - \w_{n-1}).
\end{align*}
We then associate a (generalized) Young diagram $Y$ to $\lambda$
in the following manner.  If $\w_n \ge \w_{n-1}$, let $\w_n -
\w_{n-1} = 2b_n + c_n$ with $b_n \in \Z_{\ge 0}$ and $c_n = 0$ or
$1$ and set

\begin{center}
\epsfig{file=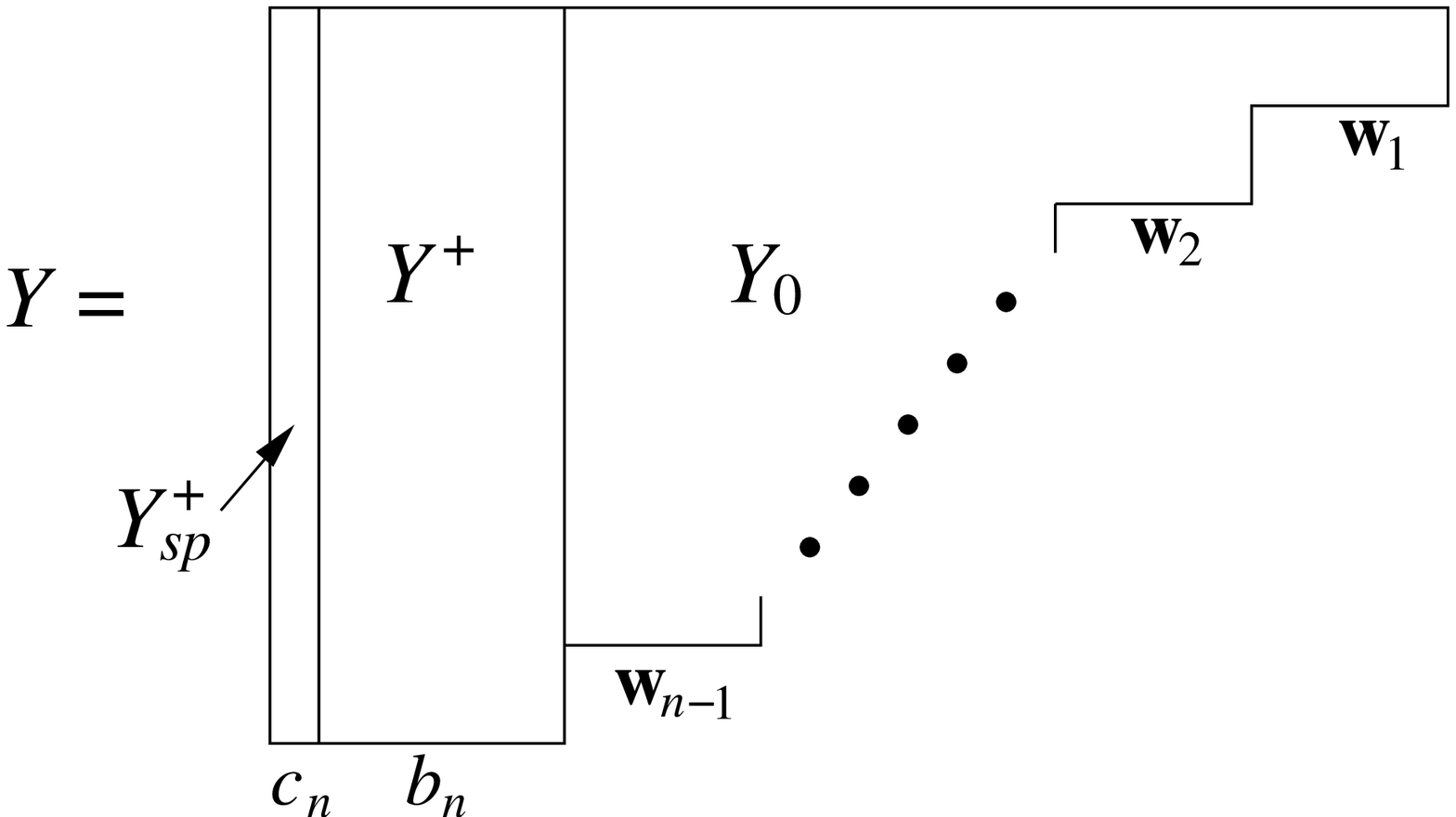,width=0.6\textwidth}
\end{center}

If $\w_n \le \w_{n-1}$, let $\w_{n-1} - \w_n = 2{\bar b}_n + {\bar
c}_n$ with ${\bar b}_n \in \Z_{\ge 0}$ and ${\bar c}_n = 0$ or $1$
and set

\begin{center}
\epsfig{file=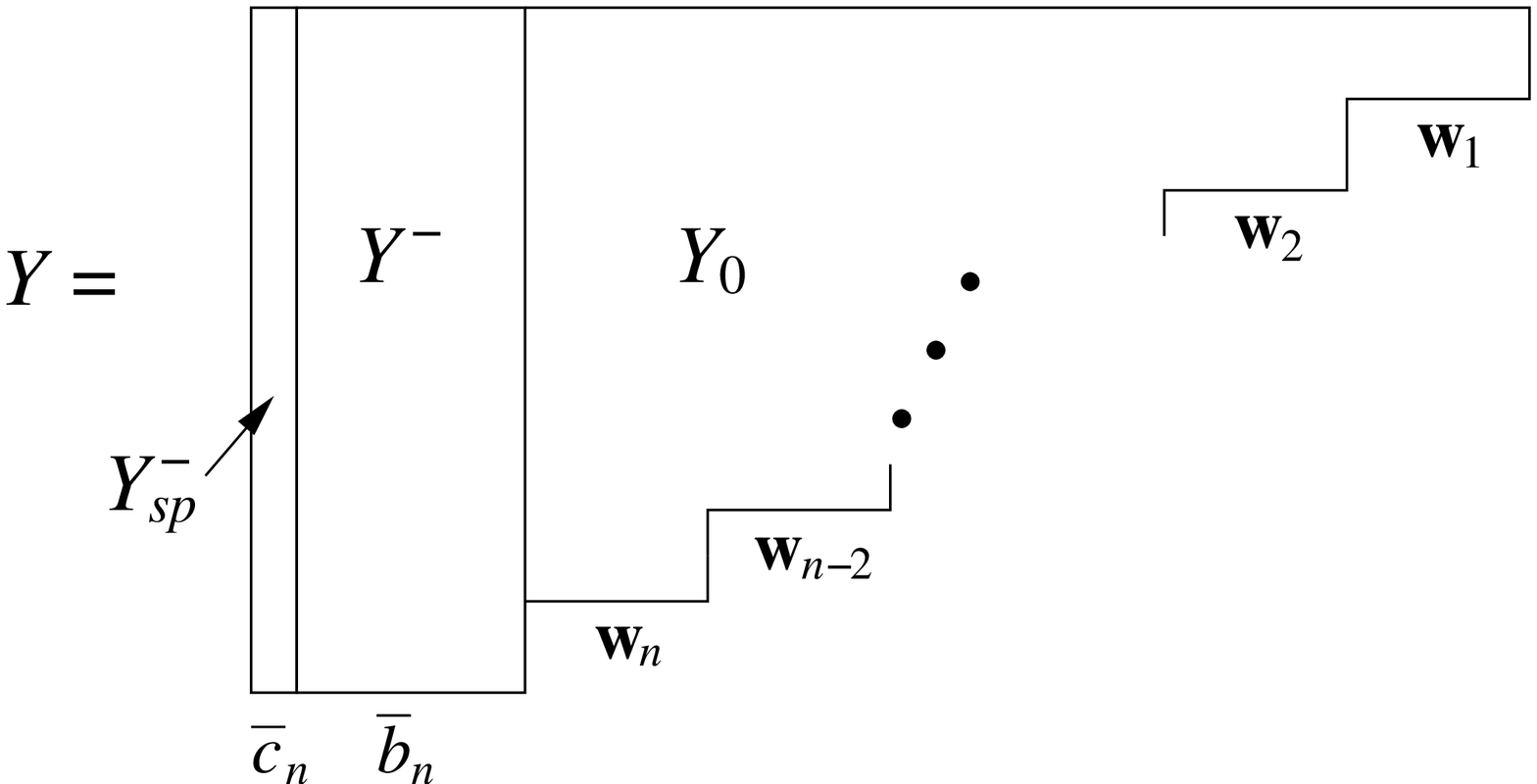,width=0.6\textwidth}
\end{center}

The column $Y_{sp}^\pm$ consists of half-boxes.  We identify the
(generalized) Young diagram with the sequence of half-integers
$Y=(\lambda_1,\dots,\lambda_n)$.  A \emph{$D_n$-tableau of shape
$Y$} is a tableau obtained from $Y$ by filling in the boxes with
entries from $\N$.  A $D_n$-tableau is said to be
\emph{semistandard} if
\begin{enumerate}
\item the entries in each row are weakly increasing and $n$ and
$\bar n$ do not appear simultaneously,
\item the entries in each
column of $Y_0$ and $Y^\pm$ are strictly increasing, but $n$ and
$\bar n$ can appear successively,
\item the entries in the column
$Y_{sp}^\pm$ are strictly increasing, and $i$ and $\bar i$ do not
appear simultaneously.
\end{enumerate}
For a $D_n$-tableau $T$, we write

\begin{center}
\epsfig{file=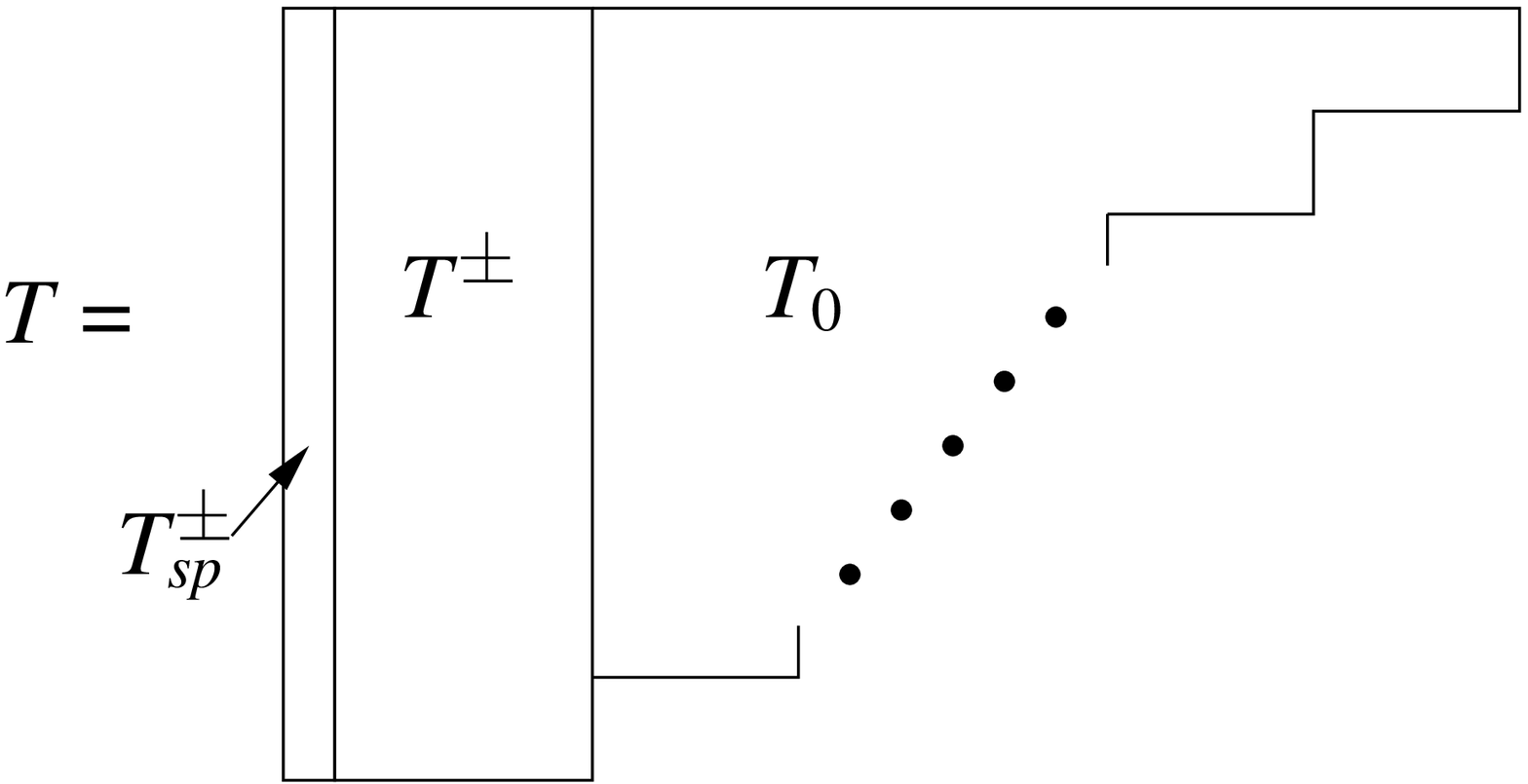,width=0.6\textwidth}
\end{center}

and we define its weight to be
\[
\wt(T) = \sum_{i=1}^n (k_i - \overline{k_i})\epsilon_i +
\frac{1}{2} \sum_{i=1}^n (l_i - \overline{l_i}) \epsilon_i,
\]
where $k_i$ (respectively $\overline{k_i}$) is the number of $i$'s
(respectively $\bar i$'s) occurring in $T_0$ and $T^\pm$, and
$l_i$ (respectively $\overline{l_i}$) is the number of $i$'s
(respectively $\bar i$'s) occurring in $T_{sp}^\pm$.  We define
$\mathcal{B}(Y)$ to be the set of all semistandard $D_n$-tableaux
$T$ of shape $Y$ satisfying conditions (D1)-(D7) of \cite{HK}.

\begin{prop}[{\cite[Thm 8.5.2]{HK}}]
For a (generalized) Young diagram $Y$ associated with a dominant
integral weight $\lambda$, $\mathcal{B}(Y)$ is isomorphic to the
crystal graph of the $U_q(\g)$-module of irreducible highest
weight $\lambda$ under the Far-Eastern reading.
\end{prop}

Note that when we compute the action of the Kashiwara operators
$\ke_i$ and $\kf_i$ for $1 \le i \le n-1$, the entries $i$ and
$\overline{i+1}$ contribute a $+$ sign and the entries $i+1$ and
$\bar i$ contribute a $-$ sign. For $\ke_n$ and $\kf_n$, the
entries $n-1$ and $n$ contribute a + sign and the entries
$\overline{n-1}$ and $\bar n$ contribute a - sign.  As for the
type $A$ case, we will say that two entries are $i$-matched if
they correspond to a $(+,-)$ pair.  Of course, we can also have a
$+$ or $-$ sign associated to the half-column $Y_{sp}^\pm$ (which
comes last in the Far-Eastern reading) and this column can thus be
$i$-matched to an entry in the rest of the tableau.

%%%%%%%%%%%%%%%%%%%%%%%%%%%%%%%%%%%%%%%%%%%%%%%%%%%%%%%%%%%%%%%%%%%%%%%

\section{Enumeration of components in type $A$}
\label{sec:A_enum}

In \cite{FS03}, the irreducible components of Nakajima's quiver
variety for type $A$ were enumerated by certain sets of Maya
diagrams. For type $A_n^{(1)}$, this enumeration matched that of a
basis for irreducible representations of the Lie algebra given in
\cite{D89}. In this section we will define for type $A_n$ a
natural 1-1 correspondence between the irreducible components of
Nakajima's quiver variety and the semistandard Young tableaux of
shape given by the highest weight.

Let $\g=\mathfrak{sl}_n$ be the Lie algebra of type $A_{n-1}$.
Then \g\ is the space of all traceless $n \times n$ matrices. The
Cartan subalgebra $\mathfrak{h}$ is spanned by the matrices
\[
H_i = e_{i,i} - e_{i+1,i+1},\ 1 \le i \le n-1
\]
where $e_{i,j}$ is the matrix with a one in entry $(i,j)$ and
zeroes everywhere else.  Thus the dual space $\mathfrak{h}^*$ is
spanned by the simple roots
\[
\alpha_i = \epsilon_i - \epsilon_{i+1},\ 1\le i \le n-1
\]
where $\epsilon_i(e_{j,j}) = \delta_{ij}$ and the fundamental
weights are given by
\[
\omega_i = \epsilon_1 + \dots + \epsilon_i,\ 1 \le i \le n-1.
\]
Consider a dominant weight $\w = \w_1 \omega_1 + \dots + \w_{n-1}
\omega_{n-1}$.  Then
\[
\w = \lambda_1 \epsilon_1 + \dots + \lambda_{n-1} \epsilon_{n-1}
\]
where $\lambda_i = \w_i + \dots + \w_{n-1}$ and so $\w$
corresponds to a partition (or Young diagram) $\lambda(\w) =
(\lambda_1 \ge \dots \ge \lambda_{n-1})$.

Let $\mathcal{T}(\lambda)$ be the set of semistandard Young
tableaux of shape $\lambda$ with entries chosen from the set
$\{1,\dots,n\}$.  Thus the entries are weakly increasing from left
to right along each row and strictly increasing down columns.  For
$T \in \mathcal{T}(\lambda)$, let $f_T(k',k)$ equal the number of
entries of $T$ in row $k'$ equal to $k+1$.  Thus, intuitively
speaking, an entry $k$ in row $k'$ corresponds to a string with
endpoints $k'$ and $k-1$ or the representation
$(\V(k',k-1),x(k',k-1))$ (where $k=k'$ yields an empty string or
the zero representation).

\[
\ffbox{k}_{k'} \longleftrightarrow \epsfig{file=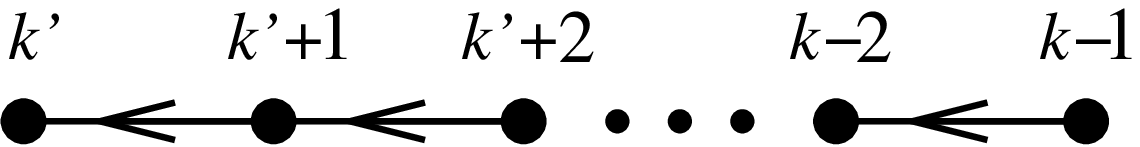,
height=.2in}
\]

\begin{prop}
The irreducible components of $\mathcal{L}(\v,\w)$ are precisely
the $X_{f_T}$ for those $f_T \in {\tilde Z}_\V$ with $T \in
\mathcal{T}(\lambda(\w))$.  Denote the component corresponding to
such a $T$ by $X_T$.  Thus, $T \leftrightarrow X_T$ is a natural
1-1 correspondence between the set $\mathcal{T}(\lambda(\w))$ and
the irreducible components of $\cup_\v \mathcal{L}(\v, \w)$.
\end{prop}

\begin{proof}
This follows from a natural generalization of Theorem~6.2 of
\cite{FS03} to the type $A_{n-1}$ case. Each of the Maya diagrams
in this theorem corresponds to a column of the Young tableau. The
representations corresponding to rows in the Maya diagram now
correspond to entries in the tableau (see
Figure~\ref{fig:maya-tab-conversion}).
\begin{figure}
\centering \epsfig{file=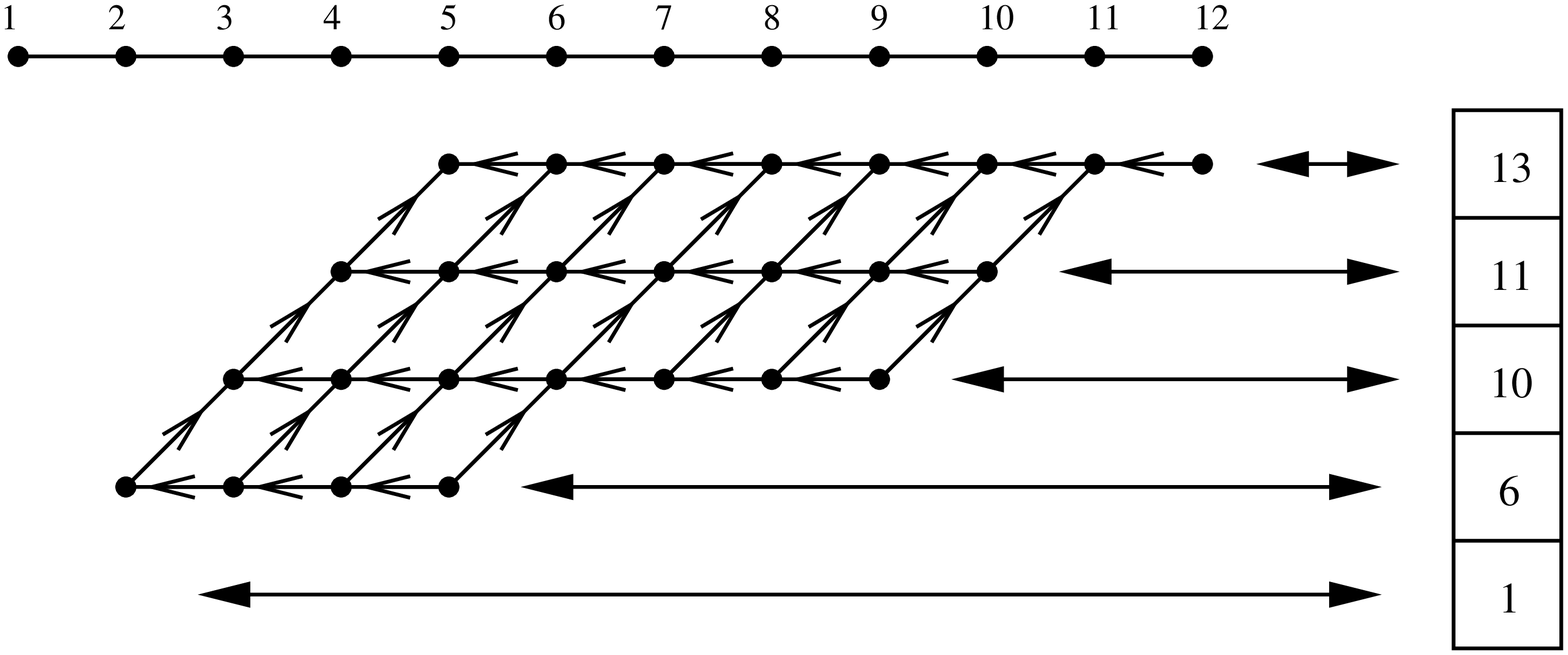, width=\textwidth}
\caption{The correspondence between strings/rows in a Maya diagram
of charge 5 and the entries in a 5 box column in a Young tableau
for the case of $A_{12}$. The entry is one greater than the degree
of the right endpoint of the string/row of the Maya diagram.  An
empty string/row corresponds to an entry equal to the where the
left endpoint of the string would have been. The column of the
Young tableau has been inverted to better show the correspondence.
\label{fig:maya-tab-conversion}}
\end{figure}
The fact that we are in the $A_{n-1}$ and not the $A_\infty$ case
results in the tableau entries being less than or equal to $n$
(since the rows of the Maya diagrams cannot extend past the
$(n-1)$st vertex).  The fact that the lengths of the rows of the
Maya diagrams weakly decrease (and hence the degrees of their
right endpoints strictly decrease) corresponds to the fact that
the entries in the tableau must strictly decrease. The fact that
the left endpoints of the strings in a Maya diagram move left by
one vertex as we move down the diagram implies that we can have at
most $k$ rows where $k$ is the charge of the Maya diagram. This
explains the fact that the corresponding column in the Young
tableau has $k$ boxes. Then the ordering on the Maya diagrams
imposed in Theorem~6.2 of \cite{FS03} corresponds to the
restriction that the entries in a tableau must weakly increase
from left to right.
\end{proof}

We now introduce some notation.  We let $\ffbox{i}_p$ denote an
element $i$ in the $p$th row of a Young tableau.  Thus, as
mentioned above, $\ffbox{i}_p$ is identified with the
representation $(\V(p,i-1), x(p,i-1))$ or the string with
endpoints $p$ and $i-1$. Recall that the case $i=p$ yields the
zero representation or the empty string.  We will often blur this
distinction and refer to the string or representation
$\ffbox{i}_p$. We call the unique vector $e_{i-1} \in \ffbox{i}_p$
of degree $i-1$ an \emph{($i-1$)-removable} vector.  Sometimes,
when we do not wish to specify the degree, we will simply say such
a vector is removable.

\begin{figure}
\centering \epsfig{file=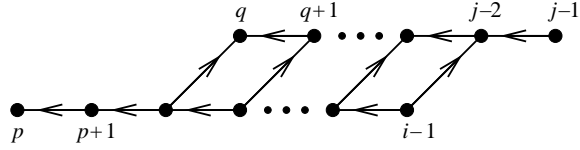, width=0.6\textwidth}
\caption{The string corresponding to an entry $i$ in the $p$th row
of a tableau can map into the string corresponding to an entry $j$
in the $q$th row if and only if $p<q \le i<j$.}
\label{fig:string_map_an}
\end{figure}

Consider a representation $x_\Omega \in \mathbf{E}_{\V,\Omega}$.
It must be a sum of indecomposable representations of the form
given in Section~\ref{sec:areps} (corresponding to strings).  Let
$(\V_1 = \V(k_1', k_1),x_1)$ and $(\V_2 = \V(k_2',k_2),x_2)$ be
two of these representations or strings. Consider the conormal
bundle to the $G_\V$-orbit through the point $x_\Omega$.  By the
proof of Theorem~5.1 of \cite{FS03}, it contains points
$(x_\Omega, x_{\bar \Omega})$ such that for some $v_1 \in \V_1$,
$x_{\bar \Omega} (v_1)$  has non-zero $v_2$-component for some
$v_2 \in \V_2$ if and only if $k_1' < k_2'$, $k_1 < k_2$, and
$k_2' \le k_1+1$. In this case we say that the string
corresponding to $(\V_1,x_1)$ maps into the string corresponding
to $(\V_2,x_2)$. In other words, the string $\ffbox{i}_p$ can map
into the string $\ffbox{j}_q$ if and only if $p<q \le i < j$.  See
Figure~\ref{fig:string_map_an}.

%%%%%%%%%%%%%%%%%%%%%%%%%%%%%%%%%%%%%%%%%%%%%%%%%%%%%%%%%%%%%%%%%%%%%

\section{Coincidence of the crystal actions: Type $A$}
\label{sec:A_ident}

As we have seen, the crystal graph of an irreducible finite
dimensional representation of $\g=\mathfrak{sl}_n$ is realized on
both the set of Young tableaux of a given shape and the set of all
irreducible components of Nakajima's quiver varieties. In this
section we show that, under the identification of tableaux with
irreducible components given in Section~\ref{sec:A_enum}, the two
crystal structures are the same.

We say that two removable vectors are $i$-matched iff their
associated entries in the tableau are $i$-matched.

\begin{lem}
\label{lem:gen-point} For a Young tableau $T \in
\mathcal{T}(\lambda(\w))$, a generic point of the irreducible
component $X_T \in B(\w)$ has a representative $(x,t)$ such that
$\ffbox{i}_p$ maps into $\ffbox{i\!+\!1}_q$ if and only if it is
$i$-matched to it.
\end{lem}

\begin{proof}
Recall that the string $\ffbox{i}_p$ can map into the string
$\ffbox{j}_q$ if and only if $p<q \le i<j$. Suppose $p_1 \le p_2 <
q \le i$ and $T$ contains entries $\ffbox{i}_{p_1}$,
$\ffbox{i}_{p_2}$ and $\ffbox{i\!+\!1}_q$.  For a generic point,
both $\ffbox{i}_{p_1}$ and $\ffbox{i}_{p_2}$ map into
$\ffbox{i\!+\!1}_q$ with non-zero coefficient.  Now,
$\ffbox{i}_{p_1}$ can map into any string $\ffbox{i}_{p_2}$ can.
Thus, by a change of basis (i.e. choosing a new representative of
the $G_\V$-orbit), we can subtract a multiple of the basis vectors
of $\ffbox{i}_{p_2}$ from those of $\ffbox{i}_{p_1}$ and assume
that in our representative $(x,t)$, $\ffbox{i}_{p_1}$ does not map
into $\ffbox{i\!+\!1}_q$. Repeating this argument, we may assume
that only the string $\ffbox{i}_p$ with maximal $p$ satisfying
$p<q$ maps into $\ffbox{i\!+\!1}_q$ (for multiple entries
$\ffbox{i}_p$ we choose the latest that is still earlier than
$\ffbox{i\!+\!1}_q$). An analogous argument then shows that we can
also assume that $\ffbox{i}_p$ only maps into the
$\ffbox{i\!+\!1}_q$ with minimal $q$ satisfying $p<q$ (for
multiple entries in the same row, we choose the earliest one that
is still later than $\ffbox{i}_p$).

Note that for two entries $\ffbox{i}_{p_1}$ and $\ffbox{i}_{p_2}$
with $p_1 < p_2$, the fact that our tableau is semistandard means
that $\ffbox{i}_{p_2}$ must occur southwest of $\ffbox{i}_{p_1}$
in $T$. Also the entry $\ffbox{i\!+\!1}_q$ for $q > p$ must either
occur (weakly) southwest of $\ffbox{i}_p$ or another entry
$\ffbox{i\!+\!1}_q$ does.  Since one entry appearing southwest of
another implies that it occurs later in the Far-Eastern reading of
the tableau, we see from the above that our chosen $\ffbox{i}_p$
and $\ffbox{i\!+\!1}_q$ are $i$-matched.  Removing them from
further consideration and repeating the above process, the result
follows.
\end{proof}

\begin{prop}
\label{prop:ec-action}
The value $\varepsilon_i(X_T)$ is equal to
the number of $i+1$ entries in $T$ which are not $i$-matched and
$\ke_i^c(X_T) = X_{T'}$ where $T'$ is obtained from $T$ by
changing all non-$i$-matched $i+1$ entries to $i$.
\end{prop}
\begin{proof}
Let $\varepsilon_i(X_T)=c$. Take a generic point $(x,t)$ of $X_T$
as in Lemma~\ref{lem:gen-point}. Then $\bigoplus_{h:\, \inc{h}=i}
x_h$ is spanned by the degree $i$ vectors of all $\ffbox{j}_p$
with $j>i+1$ and the $i$-removable vectors of those
$\ffbox{i\!+\!1}_q$ which are $i$-matched to some $\ffbox{i}_p$.
The $i$th component $\mathbf{S}_i$ of the subspace $\mathbf{S}$ in
the definition of the maps \eqref{eq:crystal-action} must be equal
to the span of these vectors and the result follows.
\end{proof}

\begin{cor}
\label{cor:fc-action} If $c$ is less than the number of
non-$i$-matched $i$ entries in $T$, then $\kf_i^c(X_T)=X_{T'}$
where $T'$ is obtained from $T$ by changing the earliest $c$
non-$i$-matched $i$ entries to $i+1$.  Otherwise $\kf_i^c(X_T)=0$.
\end{cor}
\begin{proof}
Since $\ke_i^c \kf_i^c (X_T)= X_T$ provided $\kf_i^c(X_T) \ne 0$,
we see from Proposition~\ref{prop:ec-action} that $T'$ must be
obtained from $T$ by changing $c$ of the non-$i$-matched $i$
entries to $i+1$. Furthermore, if these entries were not the
earliest then some of the new $i+1$ entries of $T'$ would be
$i$-matched and so we would not have $\ke_i^c \kf_i^c (X_T)= X_T$.
\end{proof}

\begin{theo}
\label{thm:ancrystalaction} The map from
$\mathcal{T}(\lambda(\w))$ to $B(\w)$ given by $T \mapsto X_T$ is
an isomorphism of crystals.  That is,
\begin{gather}
\label{eq:kash1} \wt(X_T)=\wt(T),\ \varepsilon_i(X_T) =
\varepsilon_i(T),\
\varphi_i(X_T) = \varphi_i(T) \\
\label{eq:kash2}
\ke_i(X_T)=X_{\ke_i(T)},\ \kf_i(X_T) = X_{\kf_i(T)}.
\end{gather}
\end{theo}
\begin{proof}
Since the equations of \eqref{eq:kash1} hold for the highest weight
element of the two crystals (the tableau with each entry equal to
its row number and the unique element of $B(\mathbf{0},\w)$), the
result will follow from the properties of a crystal if we show that
$\kf_i(X_T) = X_{\kf_i(T)}$.

Let $\varepsilon_i(X_T)=c$, $\ke_i^c(X_T)=X_{T_1}$ and
$\kf_i^{c+1}(X_{T_1})=X_{T_2}$. Note that the $i+1$ entries which
were changed to $i$ to obtain $T_1$ from $T$ must have been
earlier than those $i$ entries which are not $i$-matched.
Otherwise, these $i+1$ entries would have been $i$-matched.  Thus
it follows from Proposition~\ref{prop:ec-action} and
Corollary~\ref{cor:fc-action} that $T_2$ is obtained from $T$ by
switching the earliest non-$i$-matched $i$ entry to an $i+1$. But
this is precisely how $\kf_i(T)$ is obtained from $T$ and so
$\kf_i(T)=T_2$.
\end{proof}

%%%%%%%%%%%%%%%%%%%%%%%%%%%%%%%%%%%%%%%%%%%%%%%%%%%%%%%%%%%%%%%%%%%%%%

\section{Coincidence of the crystal actions: Type $D$}
\label{sec:D_ident}

The crystal graph of a finite dimensional representation of
$\g=\mathfrak{so}_{2n}$ is realized on Young tableaux and on the
set of irreducible components of Nakajima's quiver varieties.  In
this section we define a crystal isomorphism between the two
realizations. Our method of proof is somewhat different than that
for type $A$. For type $A$, we were able to prove independently
that the irreducible components are enumerated by the Young
tableaux appearing in the combinatorial crystal basis using the
results of \cite{FS03}. However, for type $D$, the situation is
different. While for the spin representations we have a natural
enumeration of the irreducible components by certain Young
diagrams (see \cite{S03}), which can be easily translated to the
language of tableaux as was done in Section~\ref{sec:A_enum} for
type $A$, this is not true for arbitrary integrable highest weight
representations of type $D$. The more general enumeration will
fall out from our comparison of the two presentations of the
crystal graph.

We first associate a representation of our quiver to each possible
entry in a $D_n$-tableau in $\mathcal{B}(Y)$.  To simplify our
presentation, we will describe our quiver representations with
graphs.  Each vertex represents a basis vector of the quiver
representation and will be labelled with a number 1 through $n$.
This number indicates the degree of the vector. If two vertices
are labelled $i$ and $j$ (corresponding to the basis vectors $v_i$
and $v_j$), a solid oriented edge from the vertex labelled $i$ to
the one labelled $j$ will indicate that $x_{h_{ij}}(v_i) = v_j$
and a dotted oriented edge will indicate that $x_{h_{ij}}(v_i) =
-v_j$. All other components of $x$ are zero.

We first consider the entries of $T_0$.  We make the following
identifications.

{\allowdisplaybreaks
\begin{align*}
\ffbox{i}_p &\longleftrightarrow
\epsfig{file=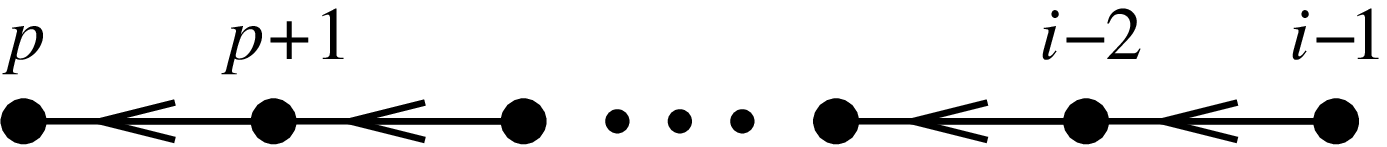,height=0.2in},\quad 1 \le p \le i \le n,\ p \ne n \\
\ffbox{\bar n}_p &\longleftrightarrow
\epsfig{file=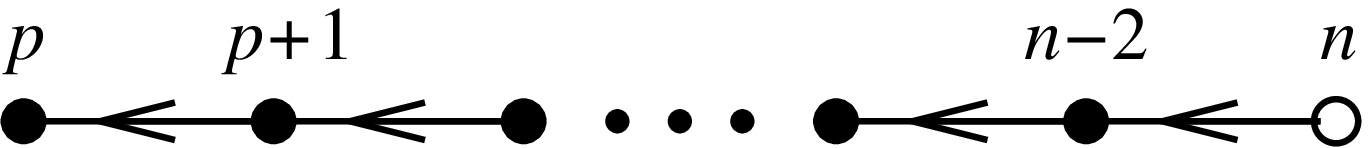,height=0.2in},\quad 1 \le p \le n-2 \\
\ffbox{\bar n}_{n-1} &\longleftrightarrow
\epsfig{file=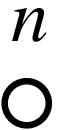,height=0.2in} \\
\ffbox{\overline{n\!-\!1}}_p &\longleftrightarrow
\raisebox{-.35in}{
\epsfig{file=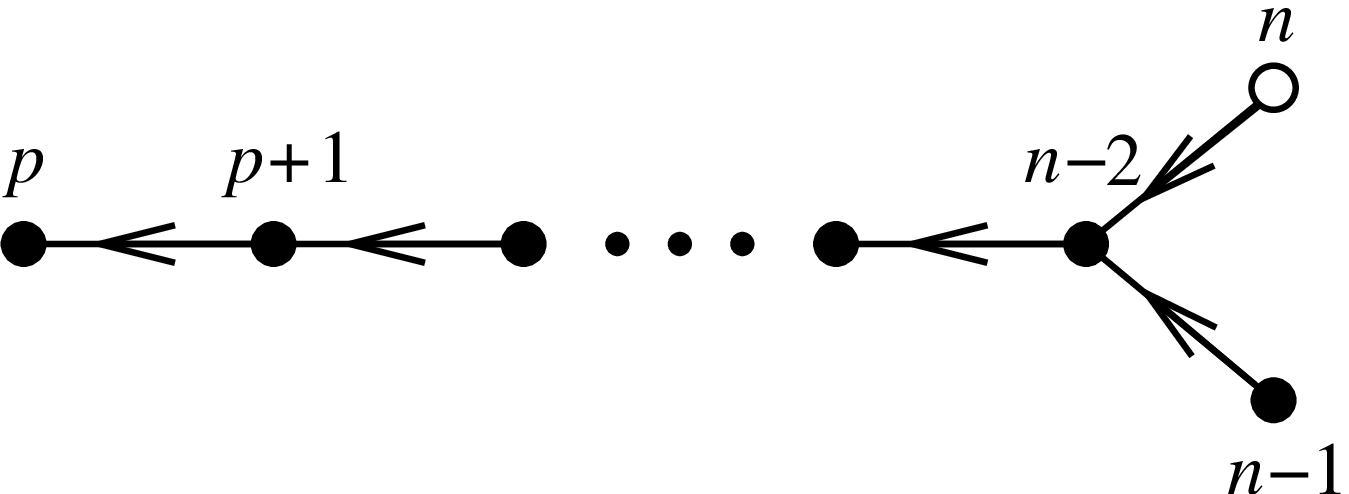,height=0.7in}},\quad 1 \le p \le n-2 \\
\ffbox{\overline{n\!-\!1}}_{n-1} &\longleftrightarrow
\raisebox{-.35in}{
\epsfig{file=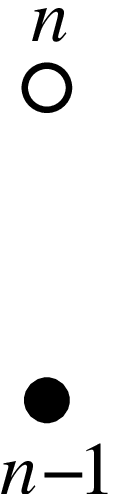,height=0.7in}} \\
\ffbox{\bar i}_p &\longleftrightarrow \raisebox{-.35in}{
\epsfig{file=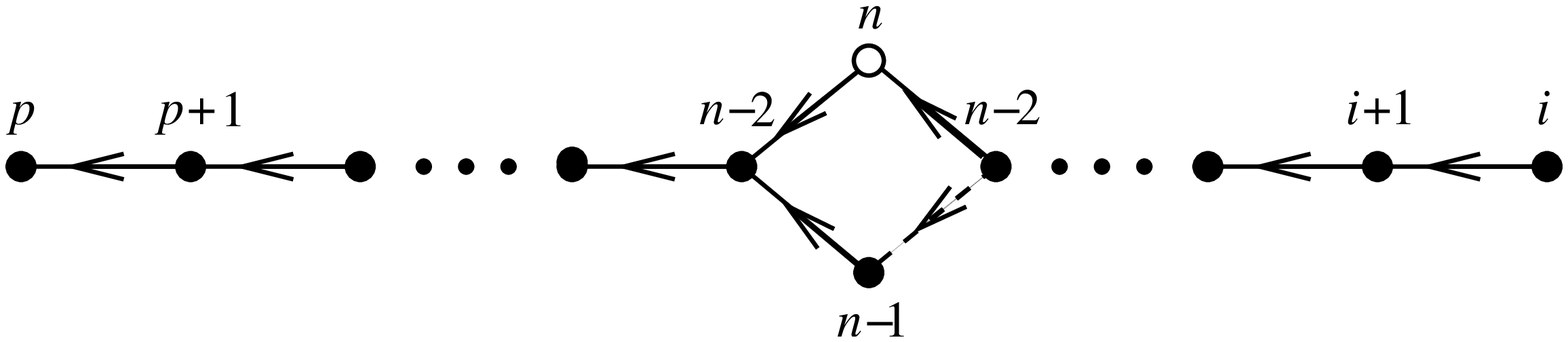,height=0.7in}},\quad 1 \le p \le n-2,\ i \le n-2 \\
\ffbox{\bar i}_{n-1} &\longleftrightarrow \raisebox{-0.35in}{
\epsfig{file=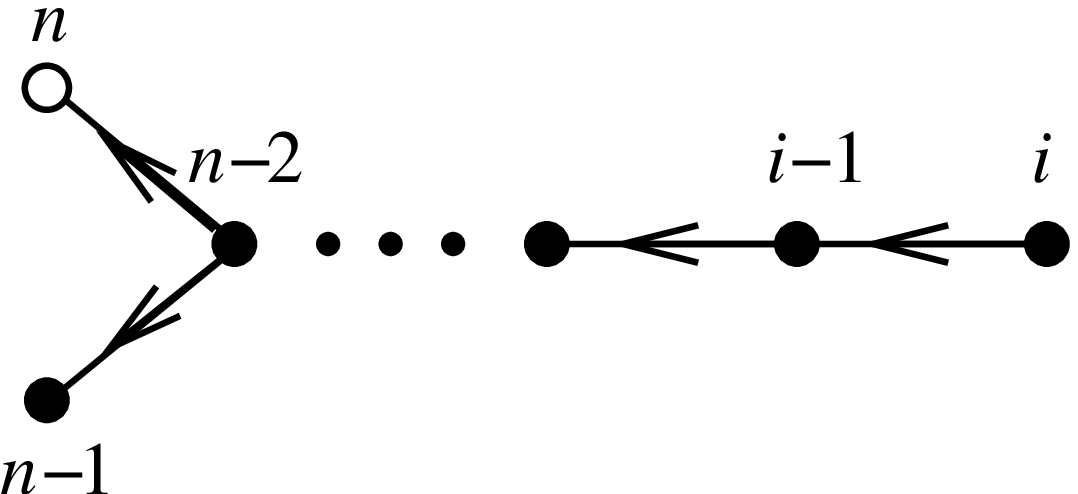,height=0.7in}},\quad i \le n-2
\end{align*}}

We use a superscript $+$ or $-$ when we want to emphasize that an
entry belongs to $T^+$/$T_{sp}^+$ or $T^-$/$T_{sp}^-$
respectively. For the entries $\ffbox{\alpha}^\pm_p$ with $1 \le p
\le n-1$, the correspondence is as for entries of $T_0$. Note that
the only possible $\ffbox{i}^+_n$ is when $i=n$ (since $i\ge n$ by
column strictness) and this corresponds to the empty string. Also,
$\ffbox{\bar n}^+_n$ is not allowed by condition (D2) of
\cite{HK}. The remaining possible entries of $T^+$ correspond to
strings as follows.

\begin{align*}
\ffbox{\bar i}_n^+ &\longleftrightarrow
\epsfig{file=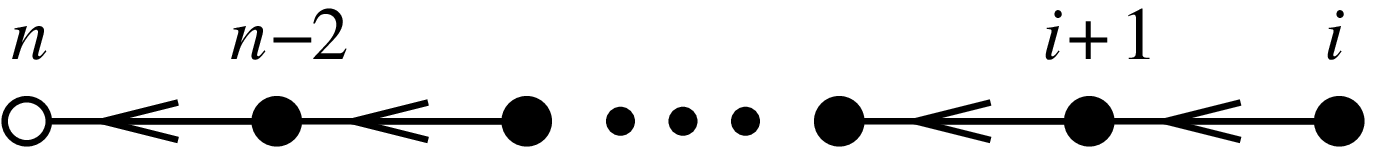,height=0.2in},\quad 1\le i \le n-2
\\
\ffbox{\overline{n\!-\!1}}_n^+ &\longleftrightarrow
\epsfig{file=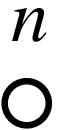,height=0.2in}
\end{align*}

Similarly, $\ffbox{i}_n^-$ is not possible for any $1 \le i \le
n-1$ by column strictness nor for $i=n$ by condition (D3) of
\cite{HK}. $\ffbox{\bar n}_n^-$ corresponds to the empty string
and the remaining possible entries of $T^-$ correspond to strings
as follows.

\[
\ffbox{\bar i}_n^- \longleftrightarrow
\epsfig{file=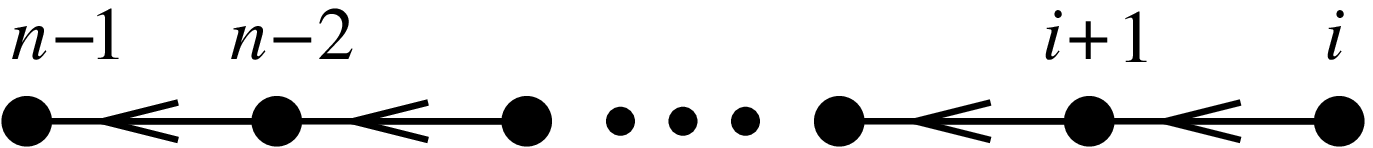,height=0.2in},\quad 1 \le i \le n-1
\]

The entries in the half boxes $\fhbox{\alpha}$ of $T_{sp}^\pm$
correspond to the empty string if $1 \le \alpha < \overline{n-1}$.
Otherwise, the entries correspond to strings as follows.

{\allowdisplaybreaks
\begin{align*}
\fhbox{\bar i}_p^+ &\longleftrightarrow
\epsfig{file=string_ibarn-.eps,height=0.2in},\quad 1 \le i \le
n-2,\ n-p \text{ is odd} \\
\fhbox{\bar i}_p^+ &\longleftrightarrow
\epsfig{file=string_ibarn+.eps,height=0.2in},\quad 1 \le i \le
n-2,\ n-p \text{ is even} \\
\fhbox{\bar i}_p^- &\longleftrightarrow
\epsfig{file=string_ibarn-.eps,height=0.2in},\quad 1 \le i \le
n-2,\ n-p \text{ is even} \\
\fhbox{\bar i}_p^- &\longleftrightarrow
\epsfig{file=string_ibarn+.eps,height=0.2in},\quad 1 \le i \le
n-2,\ n-p \text{ is odd} \\
\fsbox{\overline{n\! \text{-} \!1}}_p^+ &\longleftrightarrow
\epsfig{file=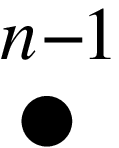,height=0.2in},\quad 1 \le i \le
n-2,\ n-p \text{ is odd} \\
\fsbox{\overline{n\! \text{-} \! 1}}_p^+ &\longleftrightarrow
\epsfig{file=string_n-1barn+.eps,height=0.2in},\quad 1 \le i \le
n-2,\ n-p \text{ is even} \\
\fsbox{\overline{n\! \text{-} \!1}}_p^- &\longleftrightarrow
\epsfig{file=string_n-1.eps,height=0.2in},\quad 1 \le i \le
n-2,\ n-p \text{ is even} \\
\fsbox{\overline{n\! \text{-} \!1}}_p^- &\longleftrightarrow
\epsfig{file=string_n-1barn+.eps,height=0.2in},\quad 1 \le i \le
n-2,\ n-p \text{ is odd}
\end{align*}}

As before, we will often blur the distinction between entries,
strings and their corresponding representations and between
vertices and the vector corresponding to them. Thus we will refer
to the strings of $T$ and say that one string is earlier than
another if its corresponding entry is, etc.

We reserve the use of the word string to refer to those
representations corresponding to a possible tableau entry. If we
can remove a vertex of degree $i$ from a string $a$ and obtain
another string, we call both the string and the vertex
\emph{$i$-removable}. When we do not wish to specify the degree of
the vertex, we may simply say that it is \emph{removable}.  It is
possible for a string to have two removable vertices.  In this
case the degree of the two vertices are necessarily $n-1$ and $n$.
Note that $\ffbox{\alpha}_p$ is $i$-removable iff it is non-empty
and $\varepsilon_i(\ffbox{\alpha})=1$. If we can add a vertex of
degree $i$ to a string $\ffbox{\alpha}_p$ and obtain another
string, we say that this string is \emph{$i$-admissible}.
Equivalently, $\ffbox{\alpha}_p$ is $i$-admissible iff
$\varphi_i(\ffbox{\alpha}) =1$.  Note that if a pair of strings is
$i$-matched then one is $i$-removable (or empty) and the other is
$i$-admissible. If a string is $i$-admissible, we call the vertex
to which the added degree $i$ vertex would map a \emph{free}
vertex. In the case that the new vertex would map into two
vertices (necessarily of degree $n-1$ and $n$), we call both
vertices free. Note that the removable and free vertex/vertices of
a string are nearly always the same.  The only case in which these
notions differ is in certain strings with a vertex of degree $n$
or $n-1$ but not both in which this vertex is the removable one
and a vertex of degree $n-2$ is the free vertex.

Note that the column of half boxes is completely determined by the
strings to which its entries correspond and that the action of the
Kashiwara operator $\kf_i$ simply corresponds to adding a vertex
of degree $i$ to a string with right endpoint of degree $i+1$
provided that no string with right endpoint of degree $i$ already
exists.  In this case we say the string is $i$-admissible.
Similarly, the action of $\ke_i$ corresponds to removing the
degree $i$ right endpoint of a string provided that no string with
right endpoint $i+1$ (or $n,n-1$ if $i=n-2$) already exits. In
this case we say the string is $i$-removable. So rather than
saying that the half column itself is $i$-matched to some entries
in the rest of the tableaux, we can say the individual strings are
$i$-matched. Because of this, we are able to treat the entries of
the half column in much the same way as the entries in the full
columns. However, when doing so, we will often treat the entries
of the half column as if they appeared in the $n$th row. This is
because of the degree of their ``terminal" vertex, which is not of
degree given by the row of the entry as in the case of the entries
of full boxes.  We also say that any string corresponding to a
half box is earlier/later than any other since in the Far-Eastern
reading, the half column is considered all at once.   For a more
detailed examination of the geometric construction of the spin
representations see \cite{S03}.

When we refer to the components of a vector in $\V$, we will
always be referring to the basis consisting of (the vectors
corresponding to) the vertices of the strings involved.  We say
that a vertex $v$ maps into a vertex $w$ if the result of applying
the map $x_h(v)$ has a non-zero $w$-component for some $h \in H$.
We also say that the string containing $v$ maps into the string
containing $w$.

\begin{lem}
\label{lem:stringmap1} If $\alpha, \beta \in \N$, $\alpha < \beta$
(where we say that $n < {\bar n}$ and ${\bar n} < n$) and $p<q$,
then $\ffbox{\alpha}_p$ can map into $\ffbox{\beta}_q$ in such a
way that if $\ffbox{\alpha}_p$ is $i$-admissible and
$\ffbox{\beta}_q$ is $i$-removable, then the free vertex/vertices
of $\ffbox{\alpha}_p$ map(s) into the $i$-removable vertex of
$\ffbox{\beta}_q$.
\end{lem}
\begin{proof}
See Figure~\ref{fig:stringmap1} for the most complicated case. The
other cases can be obtained by taking appropriate subdiagrams of
this one.
\begin{figure}
\centering \epsfig{file=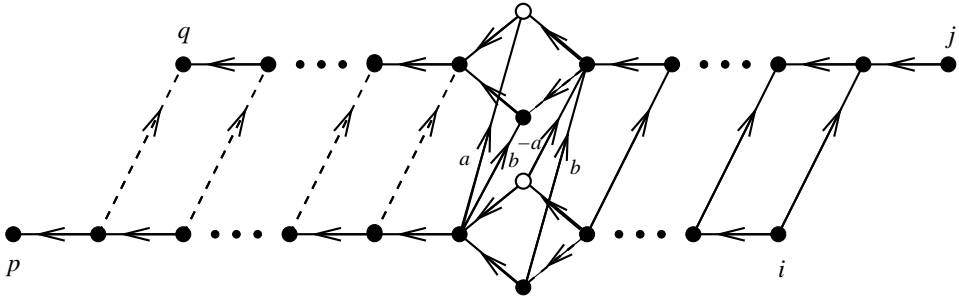,width=\textwidth}
\caption{A map from the string corresponding to an entry $\bar i$
in the $p$th row of a tableau into the string corresponding to an
entry $\bar j$ in the $q$th row for $i,j,p,q \le n-2$. Unmarked
solid and dotted lines indicate a coefficient of $1$ and $-1$
respectively. Otherwise the lines are labelled by the value of the
coefficient. These must satisfy $a+b=-1$. \label{fig:stringmap1}}
\end{figure}
\end{proof}

Lemma~\ref{lem:stringmap1} is very similar to what happened in the
type $A$ case.  However, we can also have a different type of
mapping between strings in type $D$. This can be seen in the
following proposition.

\begin{prop}
\label{prop:mappings} At least one free vertex of any
$i$-admissible string can map into the $i$-removable vertex of any
later $i$-removable string.  The strings involved may correspond
to either half or full boxes.
\end{prop}

\begin{proof}
We first deal with the case of strings corresponding to full
boxes.  Let the $i$-admissible string be $\ffbox{\alpha}_p$ and
the $i$-removable string be $\ffbox{\beta}_q$.  First suppose $1
\le \alpha, \beta < \overline{n-1}$. If $i \le n-2$ then the only
possibility is that $\alpha = i$ and $\beta = i+1$. Then since
$\ffbox{\beta}_q$ is later than $\ffbox{\alpha}_p$ and $\beta >
\alpha$, we must have $p<q$ and the result follows from
Lemma~\ref{lem:stringmap1}. If $i=n-1$ then $\alpha =n-1$ or $\bar
n$ and $\beta = n$ or $\overline{n-1}$.  In all cases,
Lemma~\ref{lem:stringmap1} again gives us the result.  If $i=n$
then $\alpha=n-1$ or $n$ and $\beta={\bar n}$ or $\overline{n-1}$.
Again, Lemma~\ref{lem:stringmap1} provides the desired result in
all cases. If $\overline{n-1} \le \alpha, \beta \le {\bar 1}$, the
result also follows from Lemma~\ref{lem:stringmap1}.

Now consider the case when $1 \le \alpha < \overline{n-1}$ and
$\overline{n-1} \le \beta \le {\bar 1}$. Then if $i \le n-2$, we
must have $\alpha=i$ and $\beta={\bar i}$. In this case
$\ffbox{\alpha}_p$ can map into $\ffbox{\beta}_q$ as in
Figure~\ref{fig:stringmap2}. If $i=n-1$, then $\alpha=n-1$ or
$\bar n$ and $\beta = \overline{n-1}$.  In both cases $\alpha <
\beta$, so $p < q$ and thus Lemma~\ref{lem:stringmap1} applies. If
$i=n$, then $\alpha=n-1$ or $n$ and $\beta = \overline{n-1}$.
Again, Lemma~\ref{lem:stringmap1} applies in either case.

\begin{figure}
\centering \epsfig{file=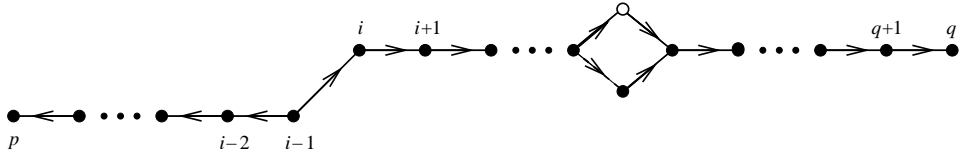,width=\textwidth}
\caption{The string corresponding to an entry $i$ in the $p$th row
of a tableau mapping into the string corresponding to an entry
$\bar i$ in the $q$th row for $i<n-2$. The second string has been
reversed from its usual presentation. \label{fig:stringmap2}}
\end{figure}

Finally, consider the case when $\overline{n-1} \le \alpha \le
{\bar 1}$ and $1 \le \beta < \overline{n-1}$. The only possibility
is $i \le n-2$. Then $\alpha=\overline{i+1}$ and $\beta=i+1$ and
$\ffbox{\alpha}_p$ can map into $\ffbox{\beta}_q$ as in
Figure~\ref{fig:stringmap3} (for $i\le n-3$) or
Figure~\ref{fig:stringmap4} (for $i=n-2$).
\begin{figure}
\centering \epsfig{file=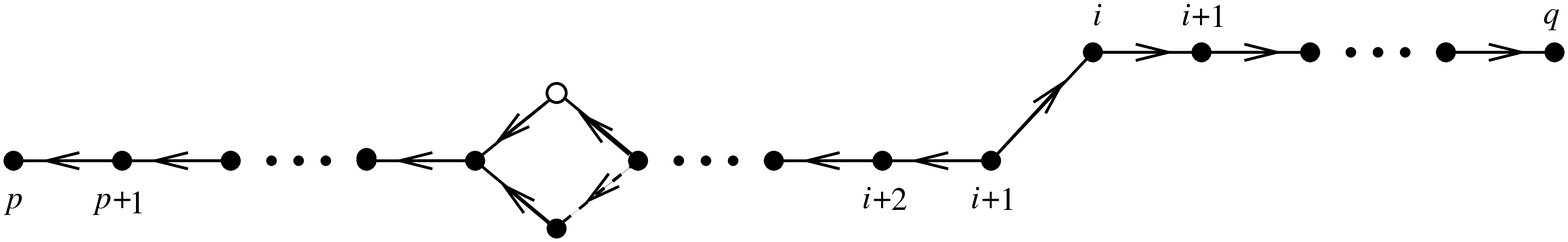,width=\textwidth}
\caption{The string corresponding to an entry $\overline{i+1}$ in
the $p$th row of a tableau mapping into the string corresponding
to an entry $i+1$ in the $q$th row for $i \le n-3$. The second
string has been reversed from its usual presentation.
\label{fig:stringmap3}}
\end{figure}
\begin{figure}
\centering \epsfig{file=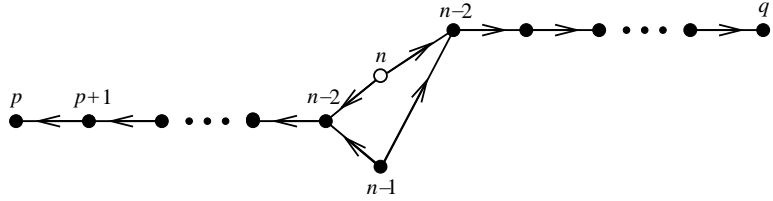,width=0.8\textwidth}
\caption{The string corresponding to an entry $\overline{n-1}$ in
the $p$th row in a tableau mapping into the string corresponding
to an entry $n-1$ in the $q$th row.  The second string has been
reversed from its usual presentation. \label{fig:stringmap4}}
\end{figure}

Now consider the case when one or both of the strings correspond
to a half box.  If both do, then the $i$-admissible string must be
$\fsbox{\overline{i\!+\!1}}^\pm$ and the $i$-removable string must
be $\fhbox{\bar i}^\pm$.  Thus the two entries must be adjacent in
the column and so one has left endpoint $n-1$ and the other $n$.
If $\fhbox{\bar i}^\pm$ has left endpoint $n$ then the map is as
in Figure~\ref{fig:halfboxmaps}.  Otherwise, we simply interchange
the labelling of the leftmost vertex of each string.
\begin{figure}
\centering \epsfig{file=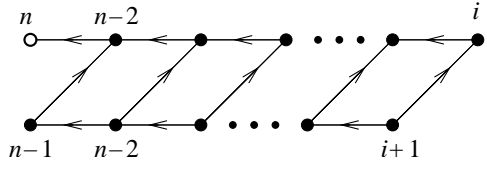,width=0.5\textwidth}
\caption{The string corresponding to the entry $\overline{i+1}$ in
a half box mapping into the string corresponding to the entry
$\bar i$ in a half box. \label{fig:halfboxmaps}}
\end{figure}

If only one of the strings corresponds to a half box, this must
necessarily be the $i$-removable string $\fhbox{i}^\pm$ since the
half entries occur later than all the full entries.  Let the
$i$-admissible string be $\ffbox{\alpha}_p$.  If $\alpha \ge
\overline{n-1}$ we can take the appropriate subdiagram of
Figure~\ref{fig:stringmap1} setting either $a$ or $b$ equal to
zero if the $i$-removable string has left endpoint $n$ or $n-1$
respectively.  If $\alpha < \overline{n-1}$ then the map can be as
in Figure~\ref{fig:fullhalfboxmap} (we picture the case where the
endpoint of the string corresponding to the half box is of degree
$n$ and $\alpha \ne {\bar n}$ -- the other cases are similar).
\begin{figure}
\centering \epsfig{file=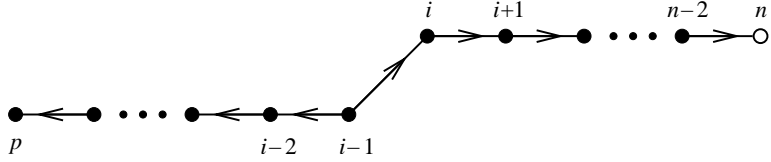,width=0.8\textwidth}
\caption{The string corresponding to the entry $i$ in the the
$p$th row of tableau mapping into the string corresponding to the
entry $\bar i$ in a half box. \label{fig:fullhalfboxmap}}
\end{figure}
\end{proof}

Now, for a given tableau $T$, let $\v^T$ be the sum of the graded
dimensions of the strings corresponding to the entries of $T$ and
let $\V^T$ be a vector space with graded dimension $\v^T$. That
is, for $i \in I$, $\v^T_i$ is equal to the number of vertices of
degree $i$ occurring in the strings of $T$.  Let $x' \in
\Lambda_{\V^T}$ be the direct sum of the representations
corresponding to the strings of $T$. We then define $A_T$ to be
the set of all $x \in \Lambda_{\V^T}$ such that for any vertex $v$
and $h \in H$, $x_h(v)$ has the same $w$-component as $x_h'(v)$
for any $w$ in the same string as $v$ and all the other components
of $x_h(v)$ lie in later strings.  Roughly speaking, we permit
strings to map only into later strings.

Let $\mathcal{C}_T$ be the union of the $G_{\V^T}$-orbits of the
points of $A_T$ and let ${\bar{\mathcal{C}}}_T$ be its closure.
Then define
\[
X_T \stackrel{\text{def}}{=} \left( \left( {\bar{\mathcal{C}}}_T
\times \sum_{i \in I} \Hom (\V^T_i, \W_i) \right) \cap
\Lambda(\v^T,\w)^{\text{st}} \right) / G_{\V^T}.
\]
Note that a priori some $X_T$ may be empty.

\begin{lem}
\label{lem:emptymatch} An $i$-admissible empty string
corresponding to a full box can only be $i$-matched to another
empty string immediately to the south of it.
\end{lem}
\begin{proof}
Unless $i=n$, an $i$-admissible empty string $a$ corresponding to
a full box must be of type $\ffbox{i}_i$.  Assume $a$ is of this
type. Suppose it is $i$-matched to another empty string $b$,
necessarily of type $\ffbox{i\!+\!1}_{i+1}$.  By column
strictness, the entry directly north of $b$ must be of type
$\ffbox{i}_i$.  Then $b$ is $i$-matched to this entry which must
therefore be the string $a$.

Now suppose $a$ is $i$-matched to the string $b$ of type
$\ffbox{\bar i}_q$.  Since $b$ is later than $a$, it is weakly
west of $a$ and $q> i$.  Then the entry in the $i$th row of the
column containing $b$ must be equal to $i$ since the tableau in
semistandard.  But then by condition (D1) of \cite{HK}, we must
have that $q$ is greater than the length of the column containing
$b$ which is a contradiction.

The only remaining possibility is that $i=n$, $a$ is of type
$\ffbox{n\!-\!1}_{n-1}$ and is $n$-matched to a string $b$ of type
$\ffbox{\bar n}_q$. Since $b$ is later than $a$ and the tableau is
semistandard, we must have $q =n$ and thus $b$ corresponds to the
empty string.

Note that the case where $i=n$, $a$ is of type $\ffbox{n}_n$ and
$b$ is of type $\ffbox{\overline{n\!-\!1}}_q$ is not possible
since we would have to have $q=n$ and then $b$ cannot be later
than $a$ by the fact that the tableau must be semistandard.
Similarly the case where $i=n-1$, $a$ is of type $\ffbox{\bar
n}_n$ and $b$ is of type $\ffbox{n}_q$ or
$\ffbox{\overline{n\!-\!1}}_q$ is not possible by the fact that
the tableau must be semistandard.
\end{proof}

\begin{prop}
\label{prop:dmatchmaps} Let $W \subset \V_i$ be the space spanned
by the $i$-removable vertices of the $i$-matched, $i$-removable
entries of a tableau $T$.  If $i=n-2$ then for each $i$-matched,
$i$-admissible string with a degree $n-2$ vertex, we extend $W$ by
the span of this vertex as well. Then a generic point of $A_T$ has
a representative such that the images under $\{x_h\ |\
\inc(h)=i\}$ of the free vertices of the $i$-admissible strings to
which these entries are $i$-matched (in the case $i=n-2$, for each
$i$-admissible string with two free vertices that does not have a
degree $n-2$ vertex we consider only the degree $n$ free vertex),
projected to $W$, form a basis of $W$.
\end{prop}

\begin{proof}
By the definition of $i$-matching, there are always more
$i$-removable strings later than an $i$-matched $i$-admissible
string than other $i$-admissible strings. Thus the result follows
from Proposition~\ref{prop:mappings}.  We need the slight
modification for the case $i=n-2$ so that the images of the free
vertices are linearly independent.
\end{proof}

\begin{cor}
\label{cor:dmatchmaps} Let $a$ be an $i$-admissible, $i$-matched
string in a tableau $T$ and take $x \in A_T$.  Suppose we extend
$x$ to a vector $v$ of degree $i$ in such a way that the resulting
point is in $\Lambda_\V$.  Then if $v$ maps into the free vertex
of $a$, either it must also map into some string earlier than $a$
or there is an $i$-admissible empty string earlier than $a$.
\end{cor}

\begin{proof}
This follows immediately from Proposition~\ref{prop:dmatchmaps}
and the moment map condition except for the case when $i=n-2$ and
one or more of the $i$-admissible, $i$-matched strings of $T$ is
of type $\ffbox{\overline{n\!-\!1}}_{n-1}$ which is the only
string with two free vertices that does not have a degree $n-2$
vertex. Let us say that $b$ is such a string and it does not have
an empty string immediately north of it (which would necessarily
be $i$-admissible). Then the string immediately north has a vector
$v'$ of degree $i=n-2$ mapping into the degree $n-1$ and $n$
vertices of $a$. Thus we can replace $v$ by a linear combination
of $v$ and $v'$ and assume that it only maps into the degree $n$
vertex (we could just have easily have used the degree $n-1$
vertices). We can do this for each such string $b$. Thus the
result again follows from Proposition~\ref{prop:dmatchmaps}.
\end{proof}

\begin{theo}
\label{thm:dncrystal}
\begin{enumerate}
\item \label{dncrystal1} The set
\[
\{X_T \ |\ T \in \mathcal{B}(Y)\}
\]
is precisely the set of irreducible components of $\cup_\v
\mathcal{L}(\v,\w)$, and

\item \label{dncrystal2} $T \mapsto X_T$ is an isomorphism of
crystals. That is,
\begin{gather}
\label{eq:dkash1} \wt(X_T)=\wt(T),\ \varepsilon_i(X_T) =
\varepsilon_i(T),\
\varphi_i(X_T) = \varphi_i(T), \\
\label{eq:dkash2} \ke_i(X_T)=X_{\ke_i(T)},\\
\label{eq:dkash3}  \kf_i(X_T) = X_{\kf_i(T)}.
\end{gather}
\end{enumerate}
\end{theo}

\begin{proof}
The variety $\mathcal{L}(0,\w)$ corresponding to the empty tableau
is a single point and equations \eqref{eq:dkash1} hold in this
case.  Also for $|\v| = \sum_i \v_i = 1$, \eqref{dncrystal1} is
true because our varieties are again points. Furthermore, it is
easy to see that \eqref{eq:dkash2} and \eqref{eq:dkash3} hold for
the highest weight tableau.  Thus, we can prove the theorem by
induction on the height $|\v|$ of $\v$.  So fix $\v$ and assume
that \eqref{dncrystal1}, \eqref{eq:dkash1}, \eqref{eq:dkash2} and
$\eqref{eq:dkash3}$ hold for all $\v'$ with $|\v'| < |\v|$.  Note
that this implies that \eqref{dncrystal1}, \eqref{eq:dkash1} and
\eqref{eq:dkash2} hold for all $\v'$ with $|\v'| \le |\v|$ by the
properties of a crystal (e.g. $\ke_i \kf_i b = b$ provided $\kf_i
b \ne 0$). Thus it is enough to show that for an irreducible
component $X_T$ of $\mathcal{L}(\v,\w)$, \eqref{eq:dkash3} holds.
Also, since we already know that both crystals are isomorphic to
that of the irreducible integrable highest weight module of
highest weight $\w$, it follows that $\kf_i(X_T)=0$ if and only if
$\kf_i(T)=0$ for any irreducible component $X_T$ of
$\mathcal{L}(\v,\w)$.  Thus we neglect this case.

By the induction hypothesis, we have that $\ke_i^c(X_T) = X_{T'}$
where $T' = \ke_i^c(T)$. Thus it suffices to show that
$\kf_i^{c+1}(X_{T'}) = X_{\kf_i^{c+1}T'}$.  Let $(\V,x)$ be a
generic point of $A_{T'}$, hence the $x$-part of a representative
of a point of $X_{T'}$. We need to describe the set of
representations $(\V',x')$ such that
\begin{gather*}
\V'_j = \V_j \ \forall\ i \ne j,\quad \V'_i = \V_i \oplus \C^{c+1}, \\
\V \text{ is $x'$-stable},\quad x'|_\V=x.
\end{gather*}
Thus we need only describe how $x'$ can act on the additional
space $\C^{c+1}$ appearing at degree $i$ in $\V'$.

Pick a basis for $\C^{c+1}$, let $v$ be the first element of the
basis, and extend $x$ generically to $v$.  So $v$ maps into all
the vertices it possibly can. We claim that if $T''=\kf_i(T')$,
then the extension of $x$ to $V \oplus \C v$ lies in $A_{T''}$.
Let $a$ be the earliest string that $v$ maps into and let it map
into the vertex $w$ of this string. Let the degree of $w$ be $j$.
Consider the following cases
\begin{enumerate}
\item \label{case1} The string $a$ has a vertex $w'$ of degree $i$
mapping into $w$.

\item \label{case2} The vertex $w$ maps into a vertex $w'$ in $a$
of degree $i$ and there is no degree $i$ vertex in $a$ mapping
into $w$.

\item \label{case3} Neither \eqref{case1} nor \eqref{case2} holds.
That is, $w$ does not map into any degree $i$ vertex in $a$ and no
degree $i$ vertex in $a$ maps into $w$.
\end{enumerate}

First consider case \eqref{case1}.  If $w'$ does not map into any
other vertex of $a$ besides $w$ then we can, by choosing a
different representative of our orbit, replace $v$ by a linear
combination of $v$ and $w'$ which does not map into $a$. So we may
assume that this case does not occur.  If $w'$ maps into another
vertex $w''$ in $a$ then we must have that the degree of $w'$ is
$n-2$ and the two vertices of $a$ it maps into are of degree $n-1$
and $n$. Then $a$ must be of type $\ffbox{\alpha}_p$ for $\alpha
\ge \overline{n-2}$ and $p \le n-1$.  Now, if $p < n-1$ then the
moment map condition implies (see
Figure~\ref{fig:dnproofdiag1.eps}) that either
\begin{figure}
\epsfig{file=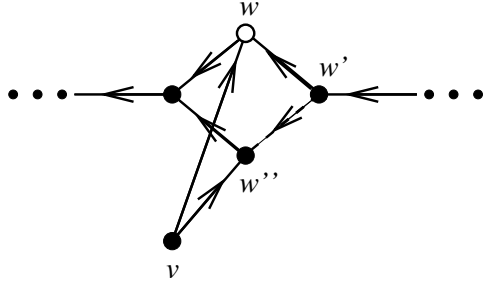, width=0.5\textwidth} \caption{The
vector $v$ mapping into the degree $n$ vertex of a string
corresponding to an entry $\alpha$, $\alpha \ge \overline{n-2}$,
in row $p$, $p<n-1$. \label{fig:dnproofdiag1.eps}}
\end{figure}
\begin{itemize}
\item $v$ maps into some vertex of a string earlier than $a$ which
in turn maps into the degree $n-2$ vertex of $a$ that $w$ maps
into, or

\item $v$ maps into $w$ and $w''$ in such a way that, considering
only $a$-components, $x_{h_{n,n-2}} x_{h_{n-2,n}} (v)
+x_{h_{n-1,n-2}} x_{h_{n-2,n-1}} (v)=0$.
\end{itemize}
The first option violates our choice of $a$ and thus we can
neglect it.  Consider the second option.  Since $w'$ satisfies the
same equation, we can again replace $v$ by a linear combination of
$v$ and $w'$ and assume that $v$ does not map into $a$.  So we can
assume that $p=n-1$. If the string directly north of $a$ is not
empty, it has a degree $n-2$ vertex $v'$ mapping into $w$ and
$w'$. We can thus replace $v$ by a linear combination of $v$, $v'$
and $w'$ (since generically the images of $v'$ and $w'$ are
linearly independent) and assume that $v$ does not map into $a$.
If the string directly north of $a$ is empty then we consider $v$
to be added to this string.

Next consider case \eqref{case2}.  We have that $x_{h_{ji}}
x_{h_{ij}} (v)$ has a non-zero $w'$-component if we consider only
$a$ and $v$.  Thus, by the moment map condition, $v$ must map into
some other vertex $w''$ mapping into $w'$. $w''$ cannot be in
another string because we picked $a$ to be the earliest string
into which $v$ maps.  Thus, $w''$ must be in the string $a$. The
only possibility is that $i=n-2$, the degree of $w$ is $n-1$ or
$n$ and the degree of $w''$ is $n$ or $n-1$ respectively. Then,
because we are assuming that there is no degree $i$ vertex mapping
into $w$, the string $a$ must be of type
$\ffbox{\overline{n\!-\!1}}$ and adding $v$ yields the string of
type $\ffbox{\overline{n\!-\!2}}$.

Finally consider case \eqref{case3}.  The only possibilities for
$w$ are that
\begin{itemize}
\item $a$ is an $i$-admissible string, $w$ is its free vertex, and
$v$ maps into it in the obvious way yielding the unique string
obtained by adding a vertex of degree $i$ to $a$,

\item $w$ is not a free vertex and is the terminal vertex of a
string (that is, $w$ maps into no other vertex in $a$), or

\item $w$ is the degree $n-2$ vertex of a string of type
$\fhbox{\bar i}^\pm$.
\end{itemize}
Consider the second possibility. Assume that $j \ne n-1,n$. It is
not possible for $a$ to be in the northernmost row in the tableau
since then $j=1$ and we must have $i=2$ (since $v$ maps into $w$).
But since $w$ is not a free vertex, there must already be a vertex
of degree $2$ mapping into it which contradicts \eqref{case3}. If
the string directly north of $a$ in the tableau is not empty then
there is another vertex $v'$ in this string of degree $i$ mapping
into $w$ in the manner of Lemma~\ref{lem:stringmap1}.  Then we
could replace $v$ by a suitable linear combination of $v$ and $v'$
and assume that $v$ does not map into $w$ after all.  So there
must be an empty string above $a$ and we consider $v$ to be added
to that string. If $j=n$ (the case $j=n-1$ is exactly analogous)
then $i=n-2$ and $a$ is either of type $\ffbox{\bar i}_{n-1}$ or
$\fhbox{\bar i}^\pm$. But since \eqref{case3} assumes that $w$ has
no vertex of degree $i$ in $a$ mapping into it and $w$ is not
free, $a$ is actually of type $\ffbox{\bar n}_{n-1}$,
$\ffbox{\overline{n\!-\!1}}_n^\pm$ or $\fsbox{\overline{n\!
\text{-} \!1}}^\pm$. If $a$ is type $\ffbox{\bar n}_{n-1}$ then
the argument used above applies and there is an empty string north
of $a$, which we consider $v$ to be added to. If $a$ is of type
$\ffbox{\overline{n\!-\!1}}_n^\pm$ or $\fsbox{\overline{n\!
\text{-} \!1}}^\pm$, then $v$ maps into the string in such a way
as to extend it to the string $\ffbox{\overline{n\!-\!2}}_n^\pm$
or $\fsbox{\overline{n\! \text{-} \!2}}^\pm$ respectively.

If $w$ is the degree $n-2$ vertex of a string of type $\fhbox{\bar
i}^\pm$, then the terminal vertex of this string must be of
different degree ($n-1$ or $n$) than $v$. By the argument used
above, there must be an $i$-admissible empty string in the half
column and we can assume that $v$ is added to this string.

Thus we have shown that earliest string that $v$ maps into (taking
the empty string in the relevant cases) is an $i$-admissible
string and it maps into the free vertex in the case the string is
non-empty.  Suppose this string is non-empty and $i$-matched. Then
by Corollary~\ref{cor:dmatchmaps}, $v$ also maps into some earlier
string. But this contradicts our choice of $a$. Thus $v$ is added
to the earliest $i$-admissible string which is not $i$-matched.
Note that this applies to empty strings by
Lemma~\ref{lem:emptymatch} since the empty strings that appeared
in the argument above were all directly north of a non-empty
string and since we can assume that $v$ is added to the earliest
non-$i$-matched empty string.

Repeating the above process for the other basis vectors of
$\C^{c+1}$, we see that after extending $x$ to all of $\C^{c+1}$,
we have changed the earliest $c+1$ non-$i$-matched $i$-admissible
entries (i.e $x' \in A_{\kf_i^{c+1}}(T')$). One should note that
for the entries of the column of half boxes, changing an entry
really corresponds to changing two entries and shuffling the order
of the entries. So we have shown that $x' \in A_{\kf_i^{c+1}T'}$.
However, to show that we obtain all of $A_{\kf_i^{c+1}T'}$, we
must show that the definition of $A_{\kf_i^{c+1}T'}$ does not
allow for any maps into $v \in \C^{c+1}$ since $V$ must be
$x'$-stable.

By definition, we need only consider strings earlier than the
string to which $v$ was added.  Suppose the vertex $w$ of an
earlier string maps into $v$. If $w$ maps into another vertex $v'$
of degree $i$ in a string $d$, then if $d$ is earlier than $a$, we
can replace $v'$ by a linear combination of $v'$ and $v$ and
assume that $w$ maps into $v'$ but not $v$. If $a$ is earlier than
$d$, we can replace $v$ by a multiple of $v$ and $v'$ and assume
that $w$ does not map into $v'$.  Continuing this process we can
assume that the only degree $i$ vertex that $w$ maps into is $v$.
But then $v$ could have mapped into $w$. Thus by our choice of the
string to which we add $v$, we must have replaced the original
(generically mapping) $v$ by a linear combination of itself and
other vectors $v_1,\dots,v_l$ so that it does not map into $w$.
Reversing this process, we can replace $v$ by a linear combination
of itself and the vectors $v_1,\dots,v_l$ so that it maps into
$w$. But then by nilpotency, $w$ cannot map into $v$. So $w$ must
have mapped into some $v_k$, $1\le k\le l$. This is a
contradiction.
\end{proof}

\begin{rem}
{\upshape It can be shown that we could have used the same
construction to define the irreducible components in finite type
$A$.  However, in type $A$ we also have the description in terms
of conormal bundles since the strings only contain edges in
$\Omega$. This description does not work in type $D$ because our
strings contain edges in both $\Omega$ and $\bar \Omega$.}
\end{rem}

\begin{rem}
{\upshape Note that in both type $A$ and $D$, our strings look
like subgraphs of the crystal graph of the vector representation.
Furthermore, the existence of certain ``maximal" maps between
strings (see Lemma~\ref{lem:stringmap1} for type $D$) gives us a
geometric interpretation of the ordering of tableaux entries.}
\end{rem}

%%%%%%%%%%%%%%%%%%%%%%%%%%%%%%%%%%%%%%%%%%%%%%%%%%%%%%%%%%%%%%%%%%%%%%%

\section{Young pyramids and the geometric crystal structure for type $A_n^{(1)}$}
\label{sec:Ahat_ident}

There exists a realization of the crystal graph of the basic
representations of the classical affine Lie algebras $A_n^{(1)}$
and $D_n^{(1)}$ (among others) in terms of combinatorial objects
known as Young walls (see \cite{HK,K03}).  It turns out that there
is a natural enumeration of the irreducible components of
Nakajima's quiver variety by Young walls.  In fact, for type
$A_n^{(1)}$, we use the geometry to suggest a generalization of
the combinatorial construction to higher level. Thus we define a
new combinatorial object which we call a \emph{Young pyramid}. On
these objects, which in level one reduce to Young walls, we will
realize the crystal of an arbitrary irreducible integrable highest
weight representation of the Lie algebra $\g =
\widehat{\mathfrak{sl}}_{n+1}$ of type $A_n^{(1)}$.

\begin{figure}
\centering \epsfig{file=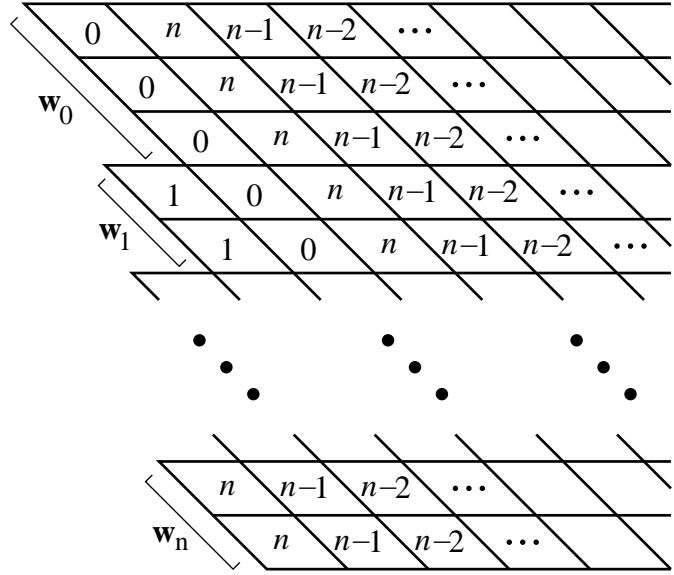,width=0.7\textwidth}
\caption{The ground state pyramid $P_\lambda$ for $\lambda =
\sum_{i=0}^n \w_i \omega_i$.  The labelling of the squares
indicates the color block that must be placed on that square.
\label{fig:groundpyramid}}
\end{figure}

Let $\lambda = \sum_{i=0}^n \w_i \omega_i$ be a dominant integral
weight.  We define $P_\lambda$ to be the \emph{ground state
pyramid} of weight $\lambda$ as shown in
Figure~\ref{fig:groundpyramid}.
We then build on the ground state by placing colored blocks (with
color $0$ through $n$).  We call a series of blocks, one on top of
the other, a \emph{stack}.  The words row and column will refer to
the rows and columns of the ground state pyramid and the blocks
placed on them.  We will use the compass directions to refer to
the relative positions of stacks (or slots in the ground state)
and the words above and below or up and down to refer to the
relative position of blocks within a stack.  The rules for
building the pyramids are as follows.
\begin{enumerate}
\item The pyramid must be built on top of the ground state
pyramid.

\item A block placed in an empty slot of the ground state wall
must be of the same color as the square on which it is placed.

\item A block placed on a block of color $i$ must be of the unique
color $j$ such that $j \equiv i+1 \mod n+1$.

\item Any stack must be weakly taller than a stack east or south
of it. That is, the height of the stacks must weakly decrease as
we move south and east.
\end{enumerate}

From now on, we assume that all colors are considered mod $n+1$.
We say a Young pyramid is \emph{proper} if the height of a stack
in the northernmost row is weakly shorter than any stack in the
southernmost row $n+1$ columns to the west. Pictorially, a Young
pyramid is proper if when we move the northernmost row to the
southern edge of the Young pyramid and shift it west by $n+1$
boxes, the condition that the heights of columns must weakly
decrease as we move south and east is still satisfied.  We say
that a Young pyramid is \emph{$n$-reduced} if for any given
height, it does not contain $n+1$ columns of that height with each
of the $n+1$ colors for the top block.  Note that in the case that
$\lambda$ is a fundamental weight, our Young pyramids are simply
the Young walls defined in \cite{HK} and all pyramids are proper.
However being $n$-reduced is different than being reduced as
defined in \cite{HK}. Let $\mathcal{F}(\lambda)$ be the set of all
proper Young pyramids built on the ground state pyramid
$P_\lambda$ and let $\mathcal{P}(\lambda)$ be the set of all
$n$-reduced proper Young pyramids in $\mathcal{F}(\lambda)$.

For two integers $k' \le k$, define $\V(k',k) \in \mathcal{V}$ to
be the vector space with basis $\{e_r\ |\ k' \le r \le k\}$.  We
require that $e_r$ has degree $i \in \{0,1,\dots,n\}$, where $r
\equiv i\, \mod n+1$.  Let $x(k',k) \in
\mathbf{E}_{\V(k',k),\Omega}$ be defined by $x(k',k) : e_r \mapsto
e_{r-1}$ for $k' \le r \le k$, where $e_{k'-1} = 0$.  Note that
the isomorphism class of this representation does not change when
$k'$ and $k$ are simultaneously translated by a multiple of $n+1$.

Let $y$ be a non-empty stack of some $P \in \mathcal{P}(\lambda)$
with bottom block of color $k'$ and height (number of boxes) $l$.
Define $k = k' + l -1$. Then let $\V^y = \V(k',k)$ and $x^y =
x(k',k)$. Equivalently, $\V^y$ is the vector space with basis
given by the blocks of the stack $y$ with the degree of each block
given by its color and $x_y$ is defined by mapping each block to
the block immediately below it (and mapping the bottom block to
zero). As for tableaux, we will use the words stack and
representation interchangeably.

Define $\V^P = \bigoplus_y \V^y$ where the sum is over all stacks
of $P$ and let $\v^P$ be its graded dimension. Then define
$x_\Omega^P \in \mathbf{E}_{\V^P,\Omega}$ to be the direct sum of
the representations $x^y$ over all the stacks of $P$. Denote by
$\mathcal{O}_P$ the $G_{\V^P}$-orbit through $x_\Omega^P$.  Let
$\mathcal{C}_P$ be the conormal bundle of $\mathcal{O}_P$ and let
${\bar{\mathcal{C}}}_P$ be its closure. Let $\W$ be a vector space
of dimension $\w$.  Then define
\[
X_P = \left( \left( {\bar{\mathcal{C}}}_P \times \sum_{j \in I}
\Hom((\V^P)_j,\W_j) \right) \cap \Lambda(\v^P,\w)^{st} \right) /
G_{\V^P}.
\]

\begin{prop}
We have that $P \leftrightarrow X_P$ is a one-to-one
correspondence between the set $\mathcal{P}(\lambda)$ and the set
of irreducible components of $\cup_\v \mathcal{L}(\v,\w)$.
\end{prop}
\begin{proof}
This follows immediately from Theorem~6.3 of \cite{FS03}. Each row
of a Young pyramid corresponds to a Maya diagram $(Y,\gamma)$. The
stacks correspond to the rows of the Young diagram $Y$ and the
charge $\gamma$ is given by the color of the slot in the ground
state pyramid of the westernmost place in the row. Then the
conditions on the heights of the columns is precisely the
condition on the ordering of the Maya diagrams and the $n$-reduced
conditions coincide.
\end{proof}

Consider a block of color $i$ in a proper Young pyramid.  The part
of the Young pyramid sitting in the same row as this block is
itself a proper Young pyramid.  We say that the block in question
is \emph{$i$-removable} if this sub-pyramid remains a proper Young
pyramid after removing this block.  Equivalently, the block is
$i$-removable if its stack is shorter than the stack immediately
to the west. Note that the entire Young pyramid may not remain a
proper Young pyramid after the removal of the block. We say that a
stack is $i$-removable if its top block is. Similarly, a place
where we may add a block of color $i$ and the sub-pyramid sitting
the same row remains a proper Young pyramid, an
\emph{$i$-admissible} slot. We say that a stack is $i$-admissible
if the top of that stack is an $i$-admissible slot.

Fix $i \in I$.  We number the stacks of $P \in
\mathcal{F}(\lambda)$, starting at zero, from tallest to shortest.
For stacks of equal height, we put $i$-removable stacks before
$i$-admissible ones (the order of the others is irrelevant) and
among $i$-removable (resp. $i$-admissible) stacks, we order them
from west to east and north to south (for stacks in the same
column).  That is, the columns to the northwest receive the lowest
indices.  When we say that one stack occurs earlier or later than
another, we are referring to this ordering. Let $y_i$ be the $i$th
stack of $P$. We assign to $y_i$ a $+$ if the stack is
$i$-admissible and a $-$ if the stack is $i$-removable. From the
sequence of $+$'s and $-$'s arranged (from left to right) in order
of increasing index, cancel out every $(+,-)$ pair to obtain a
sequence of $-$'s followed by $+$'s. This is called the
\emph{$i$-signature} of $P$. If two stacks correspond to a $(+,-)$
pair, we say they are $i$-matched. We then define $\ke_iP$ to be
the Young pyramid obtained from $P$ by removing the $i$-block from
the stack of $P$ corresponding to the rightmost $-$ in the
$i$-signature of $P$. If no $-$ exists in the $i$-signature of $P$
then $\ke_iP=0$. We define $\kf_iP$ to be the Young wall obtained
from $P$ by adding an $i$-block to the stack corresponding to the
leftmost $+$ in the $i$-signature of $P$.  If there is no $+$ in
the $i$-signature of $P$ then $\kf_iP=0$.  The astute reader will
notice that we defined an $i$-removable stack simply by
considering the sub-pyramid sitting in the row of that stack and
may worry that after removing a block from a stack $y_k$ according
to the rules above, we may not be left with a proper Young
pyramid. However, the only way this can happen is if the stack
$y_j$ immediately south of $y_k$ is of the same height and is
$i$-removable as well. But then by our numbering scheme $j > k$
and there are no $+$'s between the $-$'s corresponding to $y_k$
and $y_j$.  Thus $y_k$ is not the stack from which a block is
removed. A similar argument ensures that the procedure above for
adding a block leaves us with a proper Young pyramid as well.

Note that for the case of a fundamental weight, when the Young
pyramids are Young walls, the $i$-signature just defined is
different from the one defined in \cite{HK} because our Young wall
is reversed compared to those in \cite{HK} and thus the
cancellation of $(+,-)$ pairs occurs in the opposite direction.

Define the maps
\[
\wt : \mathcal{F}(\omega_i) \to P,\ \varepsilon_i :
\mathcal{F}(\omega_i) \to \Z,\ \varphi_i : \mathcal{F}(\omega_i)
\to \Z
\]
by
\begin{align*}
\wt(P) &= \omega_i - \sum_{j \in I} k_j \alpha_j \\
\varepsilon_i(P) &= \mbox{the number of $-$ in the $i$-signature
of $P$} \\
\varphi_i(P) &= \mbox{the number of $+$ in the $i$-signature of
$P$},
\end{align*}
where $k_i$ is the number of $i$-blocks in $P$ that have been
added to the ground-state pyramid $P_{\lambda}$.

\begin{prop}
The maps $\wt : \mathcal{F}(\lambda) \to P$, $\ke_i, \kf_i :
\mathcal{F}(\lambda) \to \mathcal{F}(\lambda) \cup \{0\}$,
$\varepsilon_i, \varphi_i : \mathcal{F}(\lambda) \to \Z$ define a
$U_q(\g)$-crystal structure on the set $\mathcal{F}(\omega_i)$ of
all proper Young pyramids.
\end{prop}
\begin{proof}
This is an straightforward verification.
\end{proof}

\begin{lem}
\label{lem:genpoint_pyramid} For a Young pyramid $P \in
\mathcal{P}(\lambda)$, a generic point of the irreducible
component $X_P \in B(\w)$ has a representative $(x,t)$ such that
\begin{enumerate}
\item The top $(i-1)$-colored block of an $i$-admissible stack
maps into the $i$-removable block of another stack if and only if
it is $i$-matched to it.

\item If there is a stack with top block of color $(i-1)$ such
that there is a stack of the same height with top block of color
$i$ immediately to the west, then this $i$-colored top block is
the only top block that the $(i-1)$-colored top block maps into.
\end{enumerate}
\end{lem}
\begin{proof}
Consider all the stacks of a given height with top block of color
$i-1$.  By the moment map condition, the $(i-1)$-colored top
blocks of these stacks can only map into the $i$-colored top
blocks of weakly shorter stacks.  If a stack with $(i-1)$-colored
top block is not admissible, it must have a stack of equal height
immediately to the west (and therefore with top block of color
$i$). By a suitable change of basis (i.e. changing our orbit
representative), we may assume that such $(i-1)$-colored top
blocks with equal height stacks immediately to the west map into
the $i$-colored top blocks of these stacks and no others.  The
stacks with $(i-1)$-colored (resp. $i$-colored) top block that do
not match up under this process are the $i$-admissible (resp.
$i$-removable) ones. Then the result follows by the argument used
in the proof of Lemma~\ref{lem:gen-point}.
\end{proof}

\begin{theo}
The map $P \mapsto X_P$ from $\mathcal{P}(\lambda)$ to the set of
irreducible components of $\cup_\v \mathcal{L}(\v,\w)$ is an
isomorphism of crystals.
\end{theo}
\begin{proof}
The proof is almost exactly analogous to that of
Theorem~\ref{thm:ancrystalaction} where we use
Lemma~\ref{lem:genpoint_pyramid} instead of
Lemma~\ref{lem:gen-point}.
\end{proof}

\begin{rem}\upshape{
Note that for the case when $\lambda$ is a fundamental weight and
our Young pyramids are Young walls, the crystal action we obtain
is different from the one in \cite{HK} because our $(+,-)$
matching occurs in the opposite direction.  However, it is also
possible to recover the crystal action of \cite{HK} from the
geometry of quiver varieties for type $A_n^{(1)}$.  One simply has
to reverse the orientation of the quiver and repeat the above
process.  Then by rotating the Young walls counterclockwise by 90
degrees, we obtain the crystal graph as presented in \cite{HK}.
However, the $n$-reduced condition is slightly more difficult to
describe in this setting and the connection to \cite{FS03} is a
little less transparent.}
\end{rem}

\begin{rem}\upshape{Note that the theory of Young pyramids could
easily be developed for type $A_\infty$ and
$\widehat{\mathfrak{gl}}_{n+1}$ using the corresponding results of
\cite{FS03} instead of the $A_n^{(1)}$ results used here.  For
both cases we would no longer have the $n$-reduced condition and
for type $A_\infty$ our set of colors would be the integers rather
than the set $\{0,1,\dots,n\}$.  The question of determining the
geometric crystal structure for $\widehat{\mathfrak{gl}}_{n+1}$
would require an extension of the results of \cite{S02} to the
varieties defined in \cite{FS03}.}
\end{rem}

%%%%%%%%%%%%%%%%%%%%%%%%%%%%%%%%%%%%%%%%%%%%%%%%%%%%%%%%%%%%%%%%%%%%%%%

\section{Young walls and the geometric crystal structure for type $D_n^{(1)}$}
\label{sec:Dhat_ident}

\begin{figure}
\centering \epsfig{file=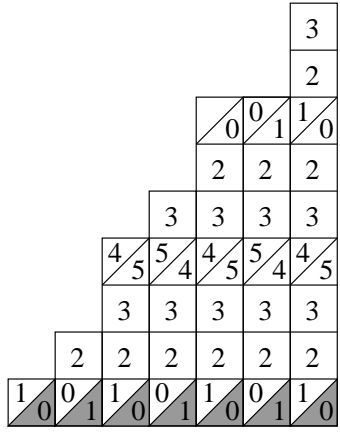,width=0.35\textwidth}
\caption{A Young wall $Y \in \mathcal{Y}(\omega_1)$ for $n=5$.
\label{fig:dnyoungwall}}
\end{figure}

Let $\g = \widehat{\mathfrak{so}}_{2n}$, $n \ge 4$, be the affine
Lie algebra of type $D_n^{(1)}$.  For a level one fundamental
weight $\omega_k$ of $\g$ (i.e. $k \in \{0,1,n-1,n\}$), define
$\mathcal{F}(\omega_k)$ as in \cite{HK} to be the set of all
proper Young walls built on the ground state wall $Y_{\omega_k}$.
Recall that a Young wall is proper if none of its full columns
have the same height.  The null root of $\g$ is given by
\[
\delta = \alpha_0 + \alpha_1 + 2\alpha_2 + \dots + 2\alpha_{n-2} +
\alpha_{n-1} + \alpha_n.
\]
A part of a column consisting of one 0-block, one 1-block, two
$i$-blocks for $2 \le i \le n-2$, one $(n-1)$-block and one
$n$-block in some cyclic order is called a \emph{$\delta$-column}.
A column in a proper Young wall $Y$ is said to contain a
\emph{removable $\delta$} if we may remove a $\delta$-column from
$Y$ and still obtain a proper Young wall.  A proper Young wall is
said to be \emph{reduced} if none of its columns contain a
removable $\delta$. Let $\mathcal{Y}(\omega_k)$ denote the set of
all reduced proper Young walls in $\mathcal{F}(\omega_k)$. An
example of an element of $\mathcal{Y}(\omega_k)$ is given in
Figure~\ref{fig:dnyoungwall}.

\begin{figure}
\centering \epsfig{file=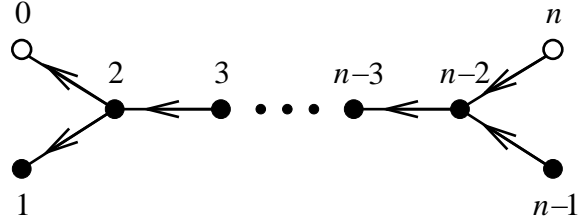,width=0.6\textwidth}
\caption{The quiver of type $D_n^{(1)}$
\label{fig:quiver_affinedn}}
\end{figure}

Let $I=\{0,1,\dots,n\}$ be the set of vertices of the Dynkin graph
of $\g$.  We label and orient the quiver as in
Figure~\ref{fig:quiver_affinedn} and we let $h_{ij}$ denote the
edge from vertex $i$ to vertex $j$.
Let $y$ be a non-empty column of some $Y \in
\mathcal{Y}(\omega_k)$. Let $\V^y$ be the vector space with basis
given by the blocks of the column $y$ with the degree of each
block given by its color and let $x^y$ be the map that sends each
box to the box immediately below it (and mapping the bottom box to
zero). In the case that a box $v$ has two boxes $w_1$ and $w_2$
below it with the degree of $w_1$ less than the degree of $w_2$,
then $v$ maps into $w_1$ with coefficient $-1$ and into $w_2$ with
coefficient $1$. Define $\V^Y = \bigoplus_y \V^y$ where the sum is
over all columns of $Y$ and let $\v^Y$ be its graded dimension.

\begin{figure}
\centering \epsfig{file=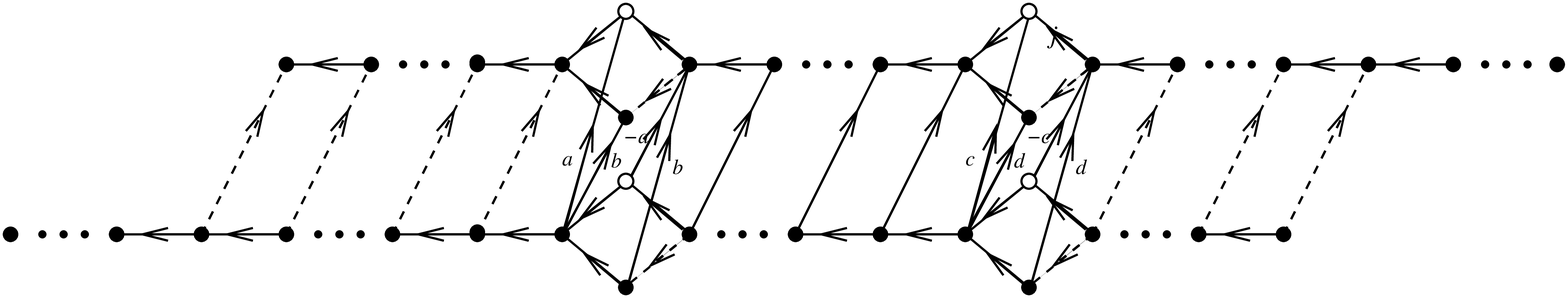,width=\textwidth}
\caption{One column mapping into another.  We depict the blocks by
vertices.  Unmarked solid and dotted lines indicate a coefficient
of $1$ and $-1$ respectively. Otherwise the lines are labelled by
the value of the coefficient. These must satisfy $a+b=-1$ and
$c+d=1$. \label{fig:affinednstringmap1}}
\end{figure}

\begin{figure}
\centering \epsfig{file=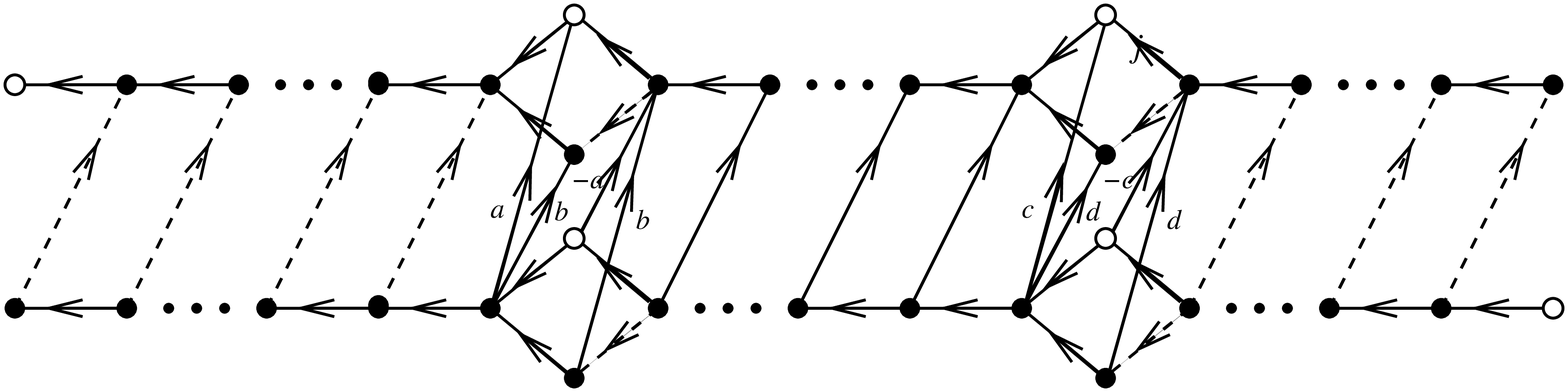,width=\textwidth}
\caption{One column mapping into another.  We depict the blocks by
vertices.  Unmarked solid and dotted lines indicate a coefficient
of $1$ and $-1$ respectively. Otherwise the lines are labelled by
the value of the coefficient. These must satisfy $a+b=-1$ and
$c+d=1$. \label{fig:affinednstringmap2}}
\end{figure}

For $Y \in \mathcal{Y}(\omega_k)$, let $x' \in \Lambda_{\V^Y}$ be
the direct sum of the representations corresponding to the columns
of $Y$.  We then define $A_Y$ to be the set of all $x \in
\Lambda_{\V^Y}$ such that for any vertex $v$ in a column $y$ and
$h \in H$, $x_h(v)$ has the same $w$-component as $x_h'(v)$ for
any $w$ in $y$ and all the other components of $x_h(v)$ lie in
columns to the right of $y$.  Roughly speaking, we permit columns
to map into other columns to their right. It is informative to
picture how one column can map into another. Two examples are
given in Figures~\ref{fig:affinednstringmap1} and
\ref{fig:affinednstringmap2}.

Let $\mathcal{C}_Y$ be the union of the $G_{\V^Y}$-orbits of the
points of $A_Y$ and let ${\bar{\mathcal{C}}}_Y$ be its closure.
Corresponding to the weight $\omega_k$, we have a dimension vector
$\w^k$.  It has $k$th component equal to one and all other
components equal to zero.  Let $\W$ be a vector space of dimension
$\w^k$ Then define
\[
X_Y \stackrel{\text{def}}{=} \left( \left( {\bar{\mathcal{C}}}_Y
\times \sum_{i \in I} \Hom (\V^Y_i, \W_i) \right) \cap
\Lambda(\v^Y,\w^k)^{\text{st}} \right) / G_{\V^Y}.
\]

We define the crystal action on the set $\mathcal{F}(\omega_k)$ as
in \cite{HK}. We call a block of color $i$ in a proper Young wall
\emph{$i$-removable} if the wall remains a proper Young wall after
removing this block. We say that a column is $i$-removable if its
top block is. Similarly, we call a place where we may add a block
of color $i$ and obtain another proper Young wall, an
\emph{$i$-admissible} slot.  We say that a column is
$i$-admissible if the top of that column is an $i$-admissible
slot.

Let $y_i$ be the $i$th column of Y, starting at zero and numbering
from right to left.  To a column $y_i$ of $Y$ we assign a $+$ if
the column is $i$-admissible and a $-$ if the column is
$i$-removable.  Combining these $+$'s and $-$'s in the same order
that the columns appear, we obtain a sequence of $+$'s and $-$'s.
From this sequence, cancel out every $(+,-)$ pair to obtain a
sequence of $-$'s followed by $+$'s.  This is called the
\emph{$i$-signature} of $Y$.  If two columns correspond to a
$(+,-)$ pair, we say they are $i$-matched. We then define $\ke_iY$
to be the Young wall obtained from $Y$ by removing the $i$-block
from the column of $Y$ corresponding to the rightmost $-$ in the
$i$-signature of $Y$. If no $-$ exists in the $i$-signature of $Y$
then $\ke_iY=0$.  We define $\kf_iY$ to be the Young wall obtained
from $Y$ by adding an $i$-block to the column corresponding to the
leftmost $+$ in the $i$-signature of $Y$.  If there is no $+$ in
the $i$-signature of $Y$ then $\kf_iY=0$.

We also define the maps
\[
\wt : \mathcal{F}(\omega_k) \to P,\ \varepsilon_i :
\mathcal{F}(\omega_k) \to \Z,\ \varphi_i : \mathcal{F}(\omega_k)
\to \Z
\]
by
\begin{align*}
\wt(Y) &= \omega_k - \sum_{j \in I} k_j \alpha_j \\
\varepsilon_i(Y) &= \mbox{the number of $-$ in the $i$-signature
of $Y$} \\
\varphi_i(Y) &= \mbox{the number of $+$ in the $i$-signature of
$Y$},
\end{align*}
where $k_i$ is the number of $i$-blocks in $Y$ that have been
added to the ground-state wall $Y_{\omega_k}$.

\begin{prop}[\cite{HK}]
The maps $\wt : \mathcal{F}(\omega_k) \to P$, $\ke_i, \kf_i :
\mathcal{F}(\omega_k) \to \mathcal{F}(\omega_k) \cup \{0\}$,
$\varepsilon_i, \varphi_i : \mathcal{F}(\omega_k) \to \Z$ define a
$U_q(\g)$-crystal structure on the set $\mathcal{F}(\omega_k)$ of
all proper Young walls.
\end{prop}

\begin{prop}[\cite{HK}]
For any $Y \in \mathcal{Y}(\omega_k)$, we have
\[
\ke_i Y \in \mathcal{Y}(\omega_k) \cup \{0\},\quad \kf_i Y \in
\mathcal{Y}(\omega_k) \cup \{0\}.
\]
Hence the set $\mathcal{Y}(\omega_k)$ has an affine crystal
structure for the quantum affine algebra $U_q(\g)$.  In fact,
$\mathcal{Y}(\omega)$ is isomorphic to the crystal of the
irreducible integrable (basic) representation of highest weight
$\omega_k$.
\end{prop}

Suppose that column $y_k$, $k>0$, is $i$-admissible and column
$y_j$ to the right of $y_k$ is $i$-removable.  Since $y_k$ is
$i$-admissible, the column $y_{k-1}$ immediately to its right must
be at least two units taller than $y_k$.  Consider the top
block(s) of $y_k$ and the block(s) of $y_{k-1}$ two rows above the
top block(s) of $y_k$. For the sake of argument, let us assume
that we are only dealing with single, full-sized blocks (the
general argument is analogous). One of these blocks is the same
color as the second block from the top in $y_j$.  Generically, the
other block maps into the top block in $y_j$.  That is,
generically, either the top block of $y_k$ or the block two rows
above in $y_{k-1}$ maps into the top block of $y_j$.  Let us call
the block that does the \emph{$y_j$-mapping} block of $y_k$.  In
the case that our block is actually composed of two half-blocks,
we call both the mapping blocks of $y_k$.  If $y_k$ is $i$-matched
to $y_j$ we simply call its $y_j$-mapping block the mapping block
(dropping the $y_j$).

If the pattern for constructing Young walls dictates that a block
of color $i$ could be added to the top of a column $y$ but the
resulting wall would not be a proper Young wall, we say that $y$
is \emph{almost $i$-admissible}.

\begin{prop}
\label{prop:genpoint_dnwall} For a Young wall $Y \in
\mathcal{Y}(\omega_k)$, a generic point of $A_Y$ has a
representative $(x,t)$ such that
\begin{enumerate}
\item \label{dn-genpoint1} The mapping block of an $i$-admissible
column maps into the $i$-removable block of the column to which it
is $i$-matched, and

\item \label{dn-genpoint2} The top block(s) of an almost
$i$-admissible column map(s) into the $i$-colored block one row
above in the column immediately to its right.
\end{enumerate}

Furthermore, let $W \subset \V_i$ be the space spanned by the
$i$-removable blocks of the $i$-matched, $i$-removable columns of
$Y$. If $i=n-2$ then for each $i$-matched, $i$-admissible column
with a degree $n-2$ block one row below the top block, we extend
$W$ by the span of this block as well. Then a generic point of
$A_Y$ has a representative such that the images under $\{x_h\ |\
\inc(h)=i\}$ of the mapping blocks of the $i$-admissible columns
to which these columns are $i$-matched (in the case $i=n-2$, for
each $i$-admissible column with two mapping blocks that does not
have a degree $n-2$ block one row below, we consider only the
degree $n$ mapping block), projected to $W$, form a basis of $W$.
\end{prop}

\begin{proof}
\eqref{dn-genpoint1} follows from the above comments and
\eqref{dn-genpoint2} follows from the existence of the type of
mappings between columns depicted in
Figures~\ref{fig:affinednstringmap1} and
\ref{fig:affinednstringmap2}.  Then the second statement follows
from the fact that by the definition of $i$-matching, there are
always more $i$-removable strings later than an $i$-matched
$i$-admissible string than other $i$-admissible strings. We need
the slight modification for the case $i=n-2$ so that the images of
the mapping blocks are linearly independent.
\end{proof}

\begin{cor}
\label{cor:genpoint_dnwall} Let $y$ be an $i$-admissible,
$i$-matched column in a Young wall $Y$ and take $x \in A_Y$.
Suppose we extend $x$ to a vector $v$ of degree $i$ in such a way
that the resulting point is in $\Lambda_\V$.  Then if $v$ maps
into the top block(s) of $y$ (or second block from the top if the
top block has half unit thickness), either it must also map into
some column to the left of $y$ or the column to the left of $y$ is
$i$-admissible and empty.
\end{cor}

\begin{proof}
This follows from Proposition~\ref{prop:genpoint_dnwall} just like
Corollary~\ref{cor:dmatchmaps} followed from
Proposition~\ref{prop:dmatchmaps}.  We use the fact that if we add
a block $v$ to a column $y$ (thus, it maps into the top block(s)
of $y$ or the second block from the top if the top of $y$ is of
half unit thickness), then since each column generically maps into
the column immediately to its right in the manner of
Figures~\ref{fig:affinednstringmap1} and
\ref{fig:affinednstringmap2}, the moment map condition implies
that $v$ must also map into the block(s) one row above it in the
column immediately to the right.  Thus $v$ maps into the mapping
block of $y$.
\end{proof}

\begin{theo}
The set $\{X_Y\ |\ Y \in \mathcal{Y}(\omega_k)\}$ is precisely the
set of irreducible components of $\cup_\v \mathcal{L}(\v,\w^k)$
and the map $Y \mapsto X_Y$ from $\mathcal{Y}(\omega_k)$ to the
set of irreducible components of $\cup_\v \mathcal{L}(\v,\w^k)$ is
an isomorphism of crystals.
\end{theo}
\begin{proof}
The proof is almost exactly analogous to that of
Theorem~\ref{thm:dncrystal} where instead of
Corollary~\ref{cor:dmatchmaps}, we use
Corollary~\ref{cor:genpoint_dnwall}. The details will be omitted.
\end{proof}

\begin{rem}
{\upshape Note that extending the above results to higher level
requires a different technique than that used in type $A_n^{(1)}$.
This is because for type $A_n^{(1)}$ every fundamental
representation is of level 1 which is not the case for type
$D_n^{(1)}$.}
\end{rem}

%%%%%%%%%%%%%%%%%%%%%%%%%%%%%%%%%%%%%%%%%%%%%%%%%%%%%%%%%%%%%%%%%%%%

\section{A connection to the path space realization}
\label{sec:paths}

In \cite{N94} and \cite{N98}, the full structure of the
representations (rather than just the crystal structure) is
defined on the homology of or a space of constructible functions
on the quiver varieties.  The irreducible components then
correspond to certain elements in these spaces, yielding a basis
for the representation.  Thus, through the identification of
irreducible components with combinatorial objects given in this
paper, we have an action of the corresponding Lie algebra \g\ on
the vector space spanned by the combinatorial objects.
Furthermore, there is a natural correspondence between Young walls
and the paths of the Kyoto model (see \cite{HK}).  One simply
reads off the top entries of columns to produce the corresponding
path.  This procedure can easily be extended to Young pyramids,
yielding a nice correspondence between Young pyramids and paths.
Then the geometric action translates into an action in the basis
given by paths.  For type $A_n^{(1)}$, an action in this basis was
given in \cite{D89b,D89} (the action for the level one case is
given explicitly while the higher level case is less explicit).
The geometric action seems to be different although the exact
connection between the two is not known.  For the case of
$D_n^{(1)}$, the definition of an action on the space of paths or
Young walls appears to be new.

%%%%%%%%%%%%%%%%%%%%%%%%%%%%%%%%%%%%%%%%%%%%%%%%%%%%%%%%%%%%%%%%%%%%%%%

\bibliographystyle{abbrv}
\bibliography{biblist}

\end{document}